\documentclass[12pt]{amsart}

\usepackage{hyperref}
\advance\textheight\topskip
\usepackage{times}

\newcommand{\figsliii}{\textup{1}}
\newcommand{\figsliv}{\textup{2}}

\usepackage[%wide,
rightcaption]{sidecap}
\usepackage{amsmath,amsfonts,amsthm,amssymb,amscd}
\usepackage[mathscr]{eucal}

\usepackage[curve,matrix,arrow]{xy}
\allowdisplaybreaks

\usepackage{times}

\newcommand{\ddiamond}{\mbox{\footnotesize$\Diamond$}}
\renewcommand{\bar}{\overline}
\renewcommand{\geq}{\geqslant}
\renewcommand{\leq}{\leqslant}

\usepackage[normalem]{ulem}
\newcommand{\opU}{\boldsymbol{\mathit{U}}}
\newcommand{\opV}{\boldsymbol{\mathit{V}}}

\newcommand{\currentU}{\mathcal{U}}

\newcommand{\nablaplus}{\nabla_{\!+}^{\vphantom{y}}}
\newcommand{\mystrut}[1]{\rule{0pt}{16pt}{\mbox{\small$#1$}}}
\newcommand{\mycirc}{\,\mbox{\footnotesize$\circ$}\,}

\newcommand{\SPsi}{\mathit{\Psi}}
\newcommand{\YXi}{\mathit{\Xi}}
\newcommand{\Tperp}{\cT_{\perp}}

\newcommand{\noo}[1]{\,{\boldsymbol{:}}#1{\boldsymbol{:}}\,}

\newcommand{\mfrac}[2]{\mbox{\small$\displaystyle\frac{#1}{#2}$}}
\newcommand{\ffrac}[2]{\raisebox{.6pt}{\mbox{\footnotesize$\displaystyle\frac{#1}{#2}$}}}
\newcommand{\fhalf}{\ffrac{1}{2}}

\newcommand{\qmodK}{\mathscr{K}}

\newcommand{\tensor}{\otimes}

\newcommand{\tq}{\widetilde{q}}
\newcommand{\Univ}{\mathscr{U}}

\newcommand{\algL}{\mathscr{L}}
\newcommand{\algN}{\mathscr{N}}
\newcommand{\algF}{\mathscr{F}}
\newcommand{\algG}{\mathscr{G}}
\newcommand{\cA}{\mathscr{A}}
\newcommand{\cR}{\mathscr{R}}
\newcommand{\polP}{\mathscr{P}}
\newcommand{\bs}[1]{\boldsymbol{#1}}
\newcommand{\wwedge}{\mathchoice{{\textstyle\bigwedge}}{\bigwedge}{\bigwedge}{\bigwedge}}
\newcommand{\ket}[1]{|#1\rangle}
\newcommand{\algA}{\mathcal{A}}
\newcommand{\algV}{\mathcal{V}}
\newcommand{\qalg}{\Univ_q\mathfrak{g}}
\newcommand{\qnilp}{\Univ_q\mathfrak{n}}
\newcommand{\cH}{\mathcal{H}}
\newcommand{\cL}{\mathcal{L}}
\newcommand{\cT}{\mathcal{T}}

\newcommand{\cW}{\mathcal{W}}
\newcommand{\cZ}{\mathcal{Z}}

\newcommand{\cX}{\mathcal{X}}
\newcommand{\cE}{\mathcal{E}}
\newcommand{\cF}{\mathcal{F}}

\providecommand{\half}{\frac{1}{2}}
\newcommand{\thalf}{\tfrac{1}{2}}
\newcommand{\dd}{\partial}
\newcommand{\ddleft}{\overleftarrow{\dd}\kern-3pt}
\newcommand{\oC}{\mathbb{C}}
\newcommand{\oN}{\mathbb{N}}

\newcommand{\oZ}{\mathbb{Z}}
\newcommand{\bref}[1]{\textbf{\ref{#1}}}
\newcommand{\ldot}{\mathbin{\boldsymbol{.}}}

\newcommand{\hSL}[1]{\widehat{s\ell}(#1)}

\newcommand{\hD}{\widehat{D}(2|1;\alpha)} 

\newcommand{\heisenberg}{\widehat{\mathfrak{h}}}
\newcommand{\hGL}[1]{\widehat{g\ell}(#1)}
\newcommand{\GL}[1]{g\ell(#1)}
\newcommand{\SL}[1]{s\ell(#1)}
\newcommand{\hSSL}[2]{\widehat{s\ell}(#1|#2)}
\newcommand{\qSSL}[2]{\Univ_q s\ell(#1|#2)}
\newcommand{\nqSSL}[2]{\mathfrak{n}_q s\ell(#1|#2)}
\newcommand{\SSL}[2]{s\ell(#1|#2)}
\newcommand{\WSL}[1]{\mathcal{W}\kern-1pts\ell(#1)}
\newcommand{\WSSL}[2]{\mathcal{W}\kern-1pts\ell(#1|#2)}

\newcommand{\BP}{W^{(2)}_3}
\newcommand{\WW}[1]{\mathscr{W}^{(2)}_{#1}}
\newcommand{\bWW}[1]{\overline{\mathscr{W}}^{(2)}_{#1}}

\numberwithin{equation}{section}
\makeatletter
\@addtoreset{equation}{section}

\def\@secnumfont{\bfseries}
\def\subsubsection{\@startsection{subsubsection}{3}%
  \z@{.5\linespacing\@plus.7\linespacing}{-.5em}%
  {\normalfont\bfseries}}
\def\paragraph{\@startsection{paragraph}{4}%
  \z@\z@{-\fontdimen2\font}%
  \normalfont\bfseries}
\def\subparagraph{\@startsection{subparagraph}{5}%
  \z@\z@{-\fontdimen2\font}%
  \normalfont\bfseries}

\makeatother

{\swapnumbers
\newtheorem{Thm}[subsection]{Theorem}

\newtheorem{lemma}[subsubsection]{Lemma}

\theoremstyle{definition}
\newtheorem{rem}[subsubsection]{Remark}
\newtheorem{example}[subsubsection]{Example} %Example

\theoremstyle{plain}

\begin{document}

\title[W2\lowercase{n} algebras]{\vspace*{-4\baselineskip}
  \mbox{}\hfill\texttt{\small\lowercase{math}.QA/\lowercase{0401164}}
  \\[3\baselineskip]
  $\WW{n}$ algebras}

\author[BL Feigin]{B.\,L.~Feigin}
  
\address{Landau Institute for Theoretical Physics, Russian Academy of
  Sciences}
  
\author[AM Semikhatov]{A.\,M.~Semikhatov}
  
\address{Theory Division, Lebedev Physics Institute, Russian Academy
  of Sciences}

\begin{abstract}
  We construct W-algebra generalizations of the $\hSL2$
  algebra\,---\,W~algebras $\WW{n}$ generated by two currents $\cE$
  and $\cF$ with the highest pole of order~$n$ in their OPE.  The
  $n\,{=}\,3$ term in this series is the Bershadsky--Polyakov $\BP$
  algebra.  We define these algebras as a centralizer (commutant) of
  the $\qSSL{n}{1}$ super quantum group and explicitly find the
  generators in a factored, ``Miura-like'' form.  Another construction
  of the $\WW{n}$ algebras is in terms of the coset
  $\hSSL{n}{1}/\hSL{n}$.  The relation between the two constructions
  involves the ``duality'' $(k\,{+}\,n\,{-}\,1)(k'\,{+}\,n\,{-}\,1)=1$
  between levels $k$ and $k'$ of two $\hSL{n}$ algebras.
\end{abstract}

\maketitle

% \setcounter{tocdepth}{2}
% \vspace*{-36pt}

% \begin{center}
%  \parbox{\textwidth}{
%    \begin{multicols}{2}
%      {\footnotesize
%        \tableofcontents}
%    \end{multicols}
%    }
% \end{center}

\thispagestyle{empty}

\section{Introduction}
The affine Lie algebra $\hSL2$ (probably the second popular in
conformal field theory after the incontestable Virasoro algebra) is
not only the first term in the sequence $\{\hSL{n}\}$ of affine Lie
algebras, but also the first term in a sequence $\{\WW{n}\}$ of W
algebras generated by two dimension-$\frac{n}{2}$ currents $\cE_n(z)$
and $\cF_n(z)$ whose operator product expansion starts with a central
term over the $n$th-order pole.
%% Although $\WW{n}$ are W, not Lie, algebras, they nevertheless
%% generalize $\hSL2$ in the sense that each $\WW{n}$ algebra is

The $\WW{n}$ algebras depend on a parameter $k\,{\in}\,\oC$.
For $n\!=\!2$, \ $k$ is the $\hSL2$ level; for $n\!=\!3$, \ $\WW{3}$
is the Bershadsky--Polyakov (BP) W~algebra $W^{(2)}_3$ and $k$ is the
level of the $\hSL3$ algebra from which it is obtained via the
``defining'' Hamiltonian reduction~\cite{[Pol],[Ber],[dBT]}.  
%% For $n\geq3$, 
Each $\WW{n}(k)$ contains a Virasoro subalgebra generated
by the energy-momentum tensor~$\cT_n(z)$ with the central charge
\begin{align}
  c_n(k)
%%   &= (2n-1)(n^2 - n -1)
%%   - n(n-1)(n-2)(k+n) - \ffrac{(n+1)n(n-1)}{k+n},
%%   \kern-8pt\\
  &=
  -\ffrac{\bigl( (k\,{+}\,n)(n \,{-}\, 1) \,{-}\, n \bigr)
    \bigl((k\,{+}\,n)(n \,{-}\, 2) n \,{-}\, n^2 \,{+}\, 1\bigr)}{
    k\,{+}\,n}
   \label{cn(k)}
\end{align}
%% (we assume $k\neq-n$ in what follows).  Each $\WW{n}(k)$ contains 
and a Heisenberg subalgebra generated by the modes of the
dimension-$1$ current $\cH_n(z)$, with the
OPEs\footnote{Normal-ordered products of composite operators are
  understood and regular terms are omitted in operator products.}
% $\cH_n(m)$, $m\,{\in}\,\oZ$, such that $[\cH_n(m),\cH_n(\ell)] =
% (k(1-\frac{1}{n}) + n-2)m\delta_{m+\ell,0}$,
\begin{gather}\label{HH-OPE}
  \begin{split}    
    \cH_n(z)\,\cH_n(w)&=
    \frac{%%\frac{1}{n}(n - 1)k + n - 2
      \ell_n(k)
    }{(z\,{-}\,w)^2},\quad
    \ell_n(k)=\ffrac{n\,{-}\,1}{n}k\,{+}\,n\,{-}\,2,\\
    \cH_n(z)\,\cE_n(w)&=\frac{\cE_n}{z\,{-}\,w},
    \qquad
    \cH_n(z)\,\cF_n(w)=-\frac{\cF_n}{z\,{-}\,w}.
  \end{split}
\end{gather}  
%% The currents $\cE_n(z)$, $\cF_n(z)$, and $\cH_n(z)$ are primary with
%% respect to $\cT_n(z)$, and 
The operator product expansion of $\cE_n$ and $\cF_n$ starts as
\begin{gather}\label{XX-short}
  \cE_n(z)\,\cF_n(w)=\frac{\lambda_{n-1}(n,k)
  }{
    (z\,{-}\,w)^n}
  {}\,{+}\,\frac{
    n
    \lambda_{n-2}(n,k)
    \,\cH_n(w)}{(z\,{-}\,w)^{n-1}}
  +\mbox{\Large$\dots$},
\end{gather}
where $\lambda_{i}(n,k)$ are numerical coefficients defined
in~\eqref{lambda-def} and the dots denote lower-order poles involving
operators of dimensions~\hbox{$\,{\geq}\,2$}.  The energy-momentum
tensor is extracted from the pole of order $n\,{-}\,2$ such that are
$\cE_{n}$ and $\cF_{n}$ are dimension-$\frac{n}{2}$ primary fields
(and $\cH$ is a dimension-$1$ primary), see~\eqref{XX-gen}.

The notation $\WW{n}$ extends the notation $W^{(2)}_3$ used for the BP
algebra, which was originally derived as a ``second'' Hamiltonian
reduction of $\hSL3$ and therefore had to be denoted similarly to but
distinctly from~$W_3$.  We tend to interpret the
superscript~${}^{(2)}$ differently, as a reminder that the
$\mathscr{W}^{(\bs{2})}_{n}$ algebras are relatives of~$\hSL{\bs{2}}$.
The lower $\WW{n}$ algebras\,---\,the Bershadsky--Polyakov algebra
$\BP$ and the otherwise not celebrated algebra~$\WW{4}$\,---\,are
described in more detail in Appendix~\ref{app:the-algebras}.

Each algebra $\WW{n}$ can be derived in at least two different ways:
from the centralizer of the quantum group $\qSSL{n}{1}$ and from the
coset theory $\hSSL{n}{1}\!\bigm/\!\hSL{n}$.

\subsubsection*{$\WW{n}$ as a centralizer of $\qSSL{n}{1}$} 
The standard definition of the W~algebra associated with a given root
system is as the centralizer of the corresponding quantum
group~\cite{[FFr]} (cf.~\cite{[BwS]}).  We recall that a general
principle in the theory of vertex-operator algebras consists in a
certain duality between vertex-operator algebras and quantum groups.
In free-field realizations of a vertex-operator algebra $\algA$, the
coresponding quantum group $\qalg$ is represented by screening
operators; more precisely, screenings are elements of~$\qnilp$, the
nilpotent subalgebra in~$\qalg$.  They can be thought of as symmetries
of the conformal field theory associated with~$\algA$: the two
algebras, $\algA$ and $\qalg$, are each other's centralizers
(commutants) in the algebra of free fields.  More precisely, the
vertex-operator algebra $\algA$ is the centralizer of $\qalg$ in the
algebra of \textit{local} operators; in turn, the centralizer of
$\algA$ in fact contains \textit{two} quantum groups, which commute
with each other and \textit{one} of which typically suffices to single
out $\algA$ as its centralizer.\footnote{\label{foot:two}For example,
  for Virasoro theories, both quantum groups are the quantum $\SL2$;
  for $\hSL2$ theories, the two quantum groups are $\qSSL21$ and
  $\Univ_{q'}\SL2$; for the coset $\hSL2\oplus\hSL2/\hSL2$, these are
  $\Univ_q D(2;1|\alpha)$ and $\Univ_{q_1^{\vphantom{y}}}\SL2
  \tensor\Univ_{q_2^{\vphantom{y}}}\SL2
  \tensor\Univ_{q_3^{\vphantom{y}}}\SL2$~\cite{[FS-D]}.}  

In the standard construction of W~algebras, the centralizer of~$\qalg$
is sought in the algebra $\algV$ of operators that are descendants of
the identity (in other words, correspond to states in the vacuum
representation).  We generalize this by first adding one extra
Heisenberg algebra (scalar current) and then \textit{extending $\algV$
  by a one-dimensional lattice vertex-operator algebra}.  We thus seek
the quantum group centralizer in $\algV_{\xi} =
\{\polP(\dd\varphi)\,e^{m\,(\xi,\varphi)}\mid m\,{\in}\,\oZ\}$, where
$\xi$ is a lightlike lattice vector, $\varphi$ is a collection of
scalar fields, $(\;{,}~)$ is the Euclidean scalar product, and $\polP$
are differential polynomials in all components of $\dd\varphi$ (and
$\dd=\frac{\dd}{\dd z}$); this involves one scalar field (Heisenberg
algebra) more than in the standard construction.  For $\xi$ chosen
specially (and for generic values of~$k$), such a centralizer of
$\qSSL{n}{1}$ is the $\WW{n}$ algebra; the required choice of the
vector $\xi$ is determined by quantum group considerations, as we
discuss in detail in the paper.

We find the $\WW{n}$ generators $(\cE_n(z),\cH_n(z),\cF_n(z))$ in the
centralizer of $\qSSL{n}{1}$ explicitly.  Before describing our
construction, we note that because $\SSL{n}{1}$ admits inequivalent
simple root systems, we have to consider different realizations
$\WW{n[m]}$, $0\,{\leq}\,m\,{\leq}\,n$, which are isomorphic
W~algebras,
% (see Sec.~\bref{sec:realizations} for the details)
but whose generators are constructed through free fields
differently.\footnote{In general, there seems to be no
  \textit{theorem} that the W~algebras obtained as centralizers of
  screenings corresponding to different root systems of the same
  algebra are isomorphic, but in all known examples, the different
  root systems lead to the same algebra.}
%% The $m=0$ and $m=n$ realizations, called \textit{maximally
%%   asymmetric}, are somewhat special; of these two, we concentrate on
%% $m=0$, because the other is obtained from the $m=0$ one by an
%% automorphism of the root system/Dynkin diagram.
%% Such realizations correspond to the different simple root systems of
%% $\SSL{n}{1}$.  
For $m=0$ and $m=n$, the respective Dynkin diagrams are (with $n$
nodes in each case)
\begin{equation*}
  \xymatrix@1{%
    {\bullet}\ar@{-}[r]
%     &{\bullet}\ar@{-}[r]
    &{\bullet}\ar@{-}[];[]+<16pt,0pt>
    &{\dots}
    &{\bullet}\ar@{-}[];[]-<16pt,0pt>
    \ar@{-}[r]
    &{\circ}
  }
\end{equation*}
and
\begin{equation*}
  \xymatrix@1{%
    {\circ}\ar@{-}[r]
    &{\bullet}\ar@{-}[];[]+<16pt,0pt>
    &{\dots}
    &{\bullet}\ar@{-}[];[]-<16pt,0pt>
    \ar@{-}[r]
    &{\bullet}
  }
\end{equation*}
which are the same.  For $1\,{\leq}\, m\,{\leq}\, n\,{-}\,1$, the corresponding
Dynkin diagrams are ($n$ nodes)
\begin{equation*}
  \xymatrix@1{%
    {\bullet}\ar@{-}[r]
%     &{\bullet}\ar@{-}[r]
    &{\bullet}\ar@{-}[];[]+<16pt,0pt>
    &{\dots}
    &{\bullet}\ar@{-}[];[]-<16pt,0pt>
    \ar@{-}[r]
    &{\circ}\ar@{-}[r]
    &{\circ}\ar@{-}[r]
    &{\bullet}\ar@{-}[];[]+<16pt,0pt>
    &{\dots}
    &{\bullet}\ar@{-}[];[]-<16pt,0pt>
  }
\end{equation*}
with $n\,{-}\,m\,{-}\,1$ black dots (even roots) to the left of the
white dots (odd roots).  To describe the generators in the centralizer
of the corresponding quantum group, we let $R^+_i$,
$i=0,\dots,n\,{-}\,m\,{-}\,1$, and $R^-_i$, $i=0,\dots,m\,{-}\,1$, be
free-field currents with the OPEs
\begin{gather}\label{RR-ope}
  \begin{aligned}
    R_i^{+}(z)\,R_i^{+}(w)&=
    R_i^{-}(z)\,R_i^{-}(w)=\ffrac{1}{(z\,{-}\,w)^2},    
    \\
    R_i^{+}(z)\,R_j^{+}(w)&=
    R_i^{-}(z)\,R_j^{-}(w)=
    \ffrac{-k\,{-}\,n\,{+}\,1}{(z\,{-}\,w)^2},\quad i\neq j,\\
    R_i^{+}(z)\,R_j^{-}(w)&=\ffrac{k\,{+}\,n\,{-}\,1}{(z\,{-}\,w)^2},
  \end{aligned}
\end{gather}
and let $Y$ be the current with the OPEs
\begin{gather}\label{RY-ope}
  R_i^{\pm}(z)\,Y(w)=\ffrac{\pm1}{(z\,{-}\,w)^2},\qquad
  Y(z)\,Y(w)=0
\end{gather}
and $\YXi$ the corresponding scalar field, such that $\dd
e^{\YXi(z)}=\noo{Ye^{\YXi}}(z)$ (with $\dd=\frac{\dd}{\dd z}$).

\begin{Thm}\label{thm:free}
  The two currents generating $\WW{n[m]}(k)$ are given by
  \begin{align*}
    \cE_{n[m]}(z)&=
    \noo{\Bigl(\!\bigl((k\,{+}\,n\,{-}\,1)\dd
      \,{+}\,R_{m-1}^-(z)
      \bigr){}\dots{}
      \bigl((k\,{+}\,n\,{-}\,1)\dd
      \,{+}\,R_{1}^-(z)\bigr)R_0^-(z)
      \Bigr)e^{\YXi(z)}},\\
    \cF_{n[m]}(z)&=
    \noo{\Bigl(\!\bigl((k\,{+}\,n\,{-}\,1)\dd
      \,{+}\,R_{n-m-1}^+(z)\bigr){}\dots{}
      \bigl((k\,{+}\,n\,{-}\,1)\dd
      \,{+}\,R_{1}^+(z)\bigr)
      R_0^+(z)\Bigr)e^{-\YXi(z)}}.
  \end{align*}
\end{Thm}
\noindent
Here, $\dd=\frac{\dd}{\dd z}$ and the action of the derivatives is
delimited by the outer brackets (i.e., the derivatives do not act on
the exponentials).  For $m=1$, only the $R_0^-$ factor is involved in
$\cE_{n[1]}$, and for $m=0$, $\cE_{n[0]}(z)=e^{\YXi(z)}$.  For
$m=n\,{-}\,1$, only the $R^+_0$ factor is involved in $\cF_{n[n-1]}$,
and for $m=n$, $\cF_{n[n]}(z)=e^{-\YXi(z)}$.

These formulas can be obtained by noting that free-field realizations
of the generators of $\WW{n[m]}(k)$ are related to free-field
realizations of the generators of $\WW{(n-1)[m]}(k+1)$.  They
generalize the following well-known situation: the bosonic
$\beta\,\gamma$ system (which is $\WW{1}$) is bosonized through two
scalar fields; this bosonization is involved in the symmetric
realization of $\hSL2$ (which is $\WW{2[1]}$), which can be obtained
by ``rebosonizing'' the Wakimoto representation.  This has an analogue
for all~$n$; in addition, relations between the maximally asymmetric
realizations of two subsequent algebras lead to the factored form of
the generators in Theorem~\bref{thm:free}.  We do not prove that the
$\cE_{n[m]}(z)$ and $\cF_{n[m]}(z)$ currents in
Theorem~\bref{thm:free} generate the entire centralizer, but we
believe that this is true for generic~$k$.

It follows from the definition that $\WW{n}$ contains the subalgebra
$\WSSL{n}{1}\tensor\Univ\heisenberg$, where $\WSSL{n}{1}$ is the
standard W~algebra associated with the $\SSL{n}{1}$ root system and
$\heisenberg$ is a Heisenberg algebra commuting with it.  The
$\WSSL{n}{1}\tensor\Univ\heisenberg$ subalgebra is merely the ``zero
momentum'' sector of $\WW{n}$, i.e., consists of elements in the
centralizer of $\qSSL{n}{1}$ in $\algV$, descendants of the identity
operator.  Introducing the generalized parafermions $\bWW{n}$ as the
quotient over the Heisenberg subalgebra, we thus have
\begin{equation*}
  \bWW{n}=\WW{n}\!\bigm/\!\heisenberg,\qquad
  \bWW{n}\supset\WSSL{n}{1}.
\end{equation*}
The algebra $\bWW{n}$ is nonlocal, but we use it because it naturally
appears in some constructions and its locality can easily be
``restored'' by tensoring with an additional free scalar field, thus
recovering $\WW{n}$.

Although this is not in the focus of our attention in this paper, we
note that $\WSSL{n}{1}$ are the ``unifying W~algebras''~\cite{[BEHHH]}
that interpolate the rank of $\WSL{\cdot}$ algebras.  The underlying
``numerology'' is as follows~\cite{[BEHHH]}.  We take the W~algebra
$\WSL{m}$, $m\,{\in}\,\oN$, and consider its $(p,p')$ minimal model.
Its central charge is given by
\begin{equation*}
  c_{p,p'}(m)=
  2 m^3\,{-}\,m\,{-}\,1\,{-}\,(m\,{-}\,1)m(m\,{+}\,1)\ffrac{p'}{p}
 \,{-}\,(m\,{-}\,1)m(m\,{+}\,1)\ffrac{p}{p'}.
\end{equation*}
The central charge of $\WW{n}(k)$ with $k = 1\,{-}\,n\,{+}\,
\ffrac{m\,{+}\,1}{n\,{-}\,1}$ satisfies the relation
\begin{equation*}
  c_n(k)\,{-}\,1 = c_{m\,{+}\,1, m\,{+}\,n}(m).
\end{equation*}
This suggests a rank--level duality of the corresponding minimal
models,
\begin{equation}\label{rank-level}
  \WSSL{n}{1}_{1 - n + \frac{m + 1}{n - 1}} = \WSL{m}_{m + n,m + 1}.
\end{equation}
The algebra of this minimal model closes on normal-ordered
differential polynomials in the currents of conformal
dimensions~$2,3,\dots,2n\,{+}\,1$~\cite{[BEHHH]}, which is often
expressed in the notation~$\mathcal{W}(2,3,\dots,2n\,{+}\,1)$ for the
model.
%% \begin{equation*}
%%   \WSSL{n}{1}_{1 - n + \frac{m + 1}{n - 1}}
%%   = \mathcal{W}(2,3,\dots,2n+1)
%% \end{equation*}
The mechanism underlying~\eqref{rank-level} was also discussed
in~\cite{[FJM]}.  The $\WSL{m}_{m + n,m + 1}$ minimal model can be
extended by taking the tensor product with a certain lattice
vertex-operator algebra generated by a free field~$\varphi$.
A~certain primary field of $\WSL{m}_{m + n,m + 1}$ can then be
``dressed'' with $e^{\alpha\varphi}$ such that the resulting operator
becomes local and can be identified with the $\cE$ current
of~$\WW{n}$; another (``dual'') primary is dressed into the $\cF$
current by $e^{-\alpha\varphi}$.  This gives ``integrable''
representations of $\WW{n}$, generalizing the integrable
representations of~$\hSL2$.

\subsubsection*{$\WW{n}$ algebras from $\protect\hSSL{n}{1}$}
Alternatively to the quantum-group description, the $\WW{n}(k)$
algebras can be given another characterization, in terms of coset
conformal field theories.  The algebras involved in the coset
construction are with the ``dual'' level $k'$ related to $k$ by
\begin{gather}\label{dual-mult}
  (k\,{+}\,n\,{-}\,1)(k'\,{+}\,n\,{-}\,1)=1.
\end{gather}
For $k\neq -n$, Eq.~\eqref{dual-mult} is equivalent to
\begin{equation}\label{dual-add}
  \mfrac{1}{k\,{+}\,n} + \mfrac{1}{k'\,{+}\,n} =1.
\end{equation}
\begin{Thm}\label{thm:coset}Let $k'$ be related to $k$
  by~\eqref{dual-mult}.  Then the coset of $\hSSL{n}1_{k'}$ with
  respect to its even subalgebra $\hGL{n}_{k'}$ is given by
%%   isomorphic to $\WW{n}(k)$
  \begin{equation}\label{the-coset}
    \hSSL{n}{1}_{k'}/\hGL{n}_{k'}
    =\WW{n[m]}(k)/\heisenberg
  \end{equation}
  (where the $\heisenberg$ subalgebra is generated by (the modes of)
  the current $\cH_n(z)$) for any $m\in[0,\dots,n]$.  In particular,
  all the $\WW{n[m]}$ algebras with $m\,{\in}\,[0,n]$ are isomorphic.
\end{Thm}
For $n=2$, we thus recover the well-known identification
$\hSL2/\heisenberg=\hSSL21/\hGL2$.

The cosets $\hSSL{n}{1}/\hSL{n}$, which are ``almost'' the left-hand
side of~\eqref{the-coset}, can be constructed in rather explicit terms
as follows.  In $\hSSL{n}{1}$, the $2n$ fermions are organized into
two $\SL{n}$ $n$-plets~$\oC^n(z)$ and~$\smash{\overset{*}{\oC}}^n(z)$.
Then $\wwedge^{\!n}\oC^n(z)$ and
$\wwedge^{\!n}\smash{\overset{*}{\oC}}^n(z)$ are in the centralizer of
the $\hSL{n}$ subalgebra.  The coset $\hSSL{n}{1}/\hSL{n}$ is
isomorphic to the centralizer of $\hSL{n}$ in $\Univ\hSSL{n}{1}$.
This is an essential simplification compared with the general case,
where a coset is usually defined as the cohomology of some BRST
operator, and only the Virasoro algebra can be realized by operators
commuting with the subalgebra.  The cosets $\hSSL{n}{1}/\hSL{n}$ is
quite close to the $\WW{n}$ algebra: after a ``correction'' by
$e^{\pm\sqrt{n}\phi(z)}$, where $\phi$ is an auxiliary scalar with the
OPE
\begin{equation}\label{phi-phi}
  \dd\phi(z)\dd\phi(w)=-\ffrac{1}{(z\,{-}\,w)^2},
\end{equation} 
$\wwedge^{\!n}\oC^n(z)$ and
$\wwedge^{\!n}\smash{\overset{*}{\oC}}^n(z)$ become the two currents
generating~$\WW{n}(k)$.

\medskip

The contents of this paper can be outlined as follows.  In
Sec.~\ref{sec:fromqg}, we study the $\WW{n}$ algebra defined as the
centralizer of $\qSSL{n}{1}$ in a certain lattice vertex-operator
algebra.  In~\bref{sec:realizations}, we define the different
\textit{realizations} of $\WW{n}$, denoted by~$\WW{n[m]}$.
In~\bref{sec:strategy}, we describe the motivation of our approach
leading to the construction of the $\WW{n[m]}$ generators in the
centralizer of the corresponding screenings.  Actual calculations,
eventually leading to Theorem~\bref{thm:free}, are given
in~\bref{sec:recursion}.  The OPEs following from
Theorem~\bref{thm:free} are considered in~\bref{sec:miura}.  In
Sec.~\ref{sec:from-affine}, we alternatively define $\WW{n}$ in terms
of the coset~$\hSSL{n}1/\hSL{n}$.  The actual statement to be proved
is in Theorem~\bref{thm:coset2} and the proof is outlined
in~\bref{sec:outline}.  It involves finding another set of screenings
representing the $\Univ_{q'}\SL{n}$ quantum group.  The quantum-group
structure is then used to construct a vertex-operator extension
of~$\WW{n}$ by the ``denominator'' $\hGL{n}$ into the ``numerator''
$\hSSL{n}{1}$, thus inverting the coset and hence showing
that~$\WW{n}$ is indeed given by the coset construction.  Some
speculations on vertex-operator extensions are given
in~\bref{sec:lifting}.

\section{$\WW{n}$ algebras from $\qSSL{n}{1}$}\label{sec:fromqg}
In this section, we define the $\WW{n}$ algebra as the centralizer of
$\qSSL{n}{1}$ in a certain lattice vertex-operator algebra.  We then
find the algebra generators in the centralizer; we do not prove that
the centralizer is thus exhausted,
%% an explicit solution for this centralizer, leading to
%% Theorem~\bref{thm:free}.  As noted above, we do not prove that this
%% solution is unique, 
but we believe that it is for generic~$k$.

\subsection{Realizations}\label{sec:realizations} To define the
$\WW{n}$ algebra as a centralizer of the screenings representing (the
nilpotent subalgebra of)~$\qSSL{n}{1}$, we proceed as follows.  We
represent $\nqSSL{n}{1}$, the nilpotent subalgebra of~$\qSSL{n}{1}$,
by operators expressed through $n+1$ scalar fields $\varphi_i$,
$i=1,\dots,n+1$.  These operators are called screenings in what
follows.  We next extend the space of differential polynomials in the
currents $\dd\varphi_1$, \dots, $\dd\varphi_{n+1}$ by adding a
one-dimensional lattice vertex-operator algebra, i.e., the operators
$e^{m(\xi,\varphi)}$, $m\,{\in}\,\oZ$, with a specially chosen
vector~$\xi$.  \textit{The $\WW{n}$ algebra is then defined as the
  centralizer of the $\nqSSL{n}{1}$-screenings in the
  space}~$\algV_{\xi}= \{\polP(\dd\varphi)\,e^{m\,(\xi,\varphi)}\mid
m\,{\in}\,\oZ\}$, where $\polP$ are differential polynomials.

To specify this in more detail, we note that the $\SSL{n}{1}$ Lie
\textit{super}algebra admits inequivalent simple root systems, and we
must therefore distinguish between centralizers of the screenings
associated with each of the inequivalent simple root systems.  To the
algebras constructed for each of the root systems, we refer as
\textit{realizations} of~$\WW{n}$, to be denoted by $\WW{n[m]}$.  We
now consider the definitions of $\WW{n[m]}$.

\subsubsection{$\boldsymbol{n{[}0{]}}$}\label{sec:n[0]} The
\textit{maximally asymmetric} realization, denoted by $\WW{n[0]}$,
corresponds to the simple root system of $\SSL{n}{1}$ represented by
the Dynkin diagram
\begin{equation*}
  \xymatrix@1{%
    {\bullet}\ar@{-}[r]
    &{\bullet}\ar@{-}[r]
    &{\bullet}\ar@{-}[];[]+<16pt,0pt>
    &{\dots}
    &{\bullet}\ar@{-}[];[]-<16pt,0pt>
    \ar@{-}[r]
    &{\circ}
  },
\end{equation*}
where filled (open) dots denote even (odd) roots.  The corresponding
% $n\times n$
Cartan matrix of $\SSL{n}{1}$ is given by
\begin{equation*}
\addtolength{\arraycolsep}{-3pt}
  \mbox{\small$\displaystyle\begin{pmatrix}
    2&-1&0&\hdotsfor{3}&0\\
    -1&2&-1&0&\hdotsfor{2}&0\\
    \hdotsfor{7}\\
    \hdotsfor{7}\\
    0&\hdotsfor{2}&0&-1&2&-1\\
    0&\hdotsfor{3}&0&-1&0
  \end{pmatrix}.$}
\end{equation*}
{}From this Cartan matrix, we construct the screenings representing
$\nqSSL{n}{1}$.  For this, we first introduce vectors $\bs{a}_{n-1}$,
\dots, $\bs{a}_1$, $\bs{\psi}$, in~$\oC^n$ whose Gram matrix (pairwise
scalar products) is given by the ``dressed Cartan matrix''
\begin{gather}\label{Cartan-asym-dressed}
  \mbox{\addtolength{\arraycolsep}{-1pt}
    \small$\displaystyle\begin{pmatrix}
    2(k\,{+}\,n)
    \rlap{\kern-100pt$\smash{\begin{array}[t]{r}
          \bs{a}_{n-1}\\\bs{a}_{n-2}\\ \\ \\\bs{a}_{1}\\\bs{\psi}
        \end{array}\rule[-80pt]{.5pt}{86pt}}$}
    &-k\,{-}\,n&0&\hdotsfor{3}&0\\
    -k\,{-}\,n&2(k\,{+}\,n)&-k\,{-}\,n&0&\hdotsfor{2}&0\\
    \hdotsfor{7}\\
    \hdotsfor{7}\\
    0&\hdotsfor{2}&0&-k\,{-}\,n&2(k\,{+}\,n)&-k\,{-}\,n\\
    0&\hdotsfor{3}&0&-k\,{-}\,n&1
  \end{pmatrix}.$}\kern-60pt
\end{gather}
where the leftmost column indicates labeling of rows.  The determinant
of this matrix is~$-(k\,{+}\,n)^{n-1}n\ell_n(k)$, and the vectors are
therefore determined uniquely modulo a common rotation for
$k\,{\in}\,\oC\setminus\{-n,-\frac{n(n-2)}{n-1}\}$.

Let $\bs{\varphi}$ be the $n$-plet of scalar fields with the OPEs
\begin{gather*}
  \dd\varphi_i(z)\,\dd\varphi_j(w)=\ffrac{\delta_{i,j}}{(z\,{-}\,w)^2}.
\end{gather*}
With the above vectors $\bs{a}_{n-1}$, \dots, $\bs{a}_1$, $\bs{\psi}$,
we define the operators
\begin{gather*}
  E_i=\oint e^{\bs{a}_i\ldot\bs{\varphi}},~i=1,\dots,n\,{-}\,1,
  \qquad
  \SPsi=\oint e^{\bs{\psi}\ldot\bs{\varphi}},
\end{gather*}
where 
%% $\bs{\varphi}$ is an $n$-component scalar field and 
the dot denotes the Euclidean scalar product in~$\oC^n$.  These
\textit{screening operators} represent $\nqSSL{n}{1}$.  The screenings
$E_i$ are said to be \textit{bosonic} and $\SPsi$ \textit{fermionic}.

We next define the vector $\bs{\xi}\,{\in}\,\oC^n$ by its scalar
products with $\bs{a}_1$, \dots, $\bs{a}_{n-1}$, $\bs{\psi}$,
\begin{equation}\label{xi-cond-0}
  \begin{split}
    \bs{\xi}\ldot\bs{a}_i &=0,\quad i=1,\dots,n\,{-}\,1,\\
    \bs{\xi}\ldot\bs{\psi}&=1
  \end{split}
\end{equation}
and set $\bar{\algV}_{\bs{\xi}} =
\{\polP(\dd\bs{\varphi})\,e^{m\,\bs{\xi}\ldot\bs{\varphi}}\mid
m\,{\in}\,\oZ\}$.  By definition, the $n[0]$ realization of~$\bWW{n}$,
denoted by $\bWW{n[0]}$, is the centralizer
of~$(E_i)_{i=1,\dots,n\,{-}\,1}$ and~$\SPsi$
in~$\bar{\algV}_{\bs{\xi}}$.  Restoring locality\,---\,i.e.,
reconstructing $\WW{n[0]}$\,---\,requires one scalar field more.

To construct $\WW{n}$, we embed $\oC^n$ in $\oC^{n+1}$ as a coordinate
hyperplane, let $\psi$, $a_1$, \dots, $a_{n-1}$ denote the respective
images of $\bs{\psi}$, $\bs{a}_1$, \dots, $\bs{a}_{n-1}$, and extend
the set of free scalar fields accordingly, to
$\dd\varphi=\{\dd\bs{\varphi},\dd\varphi_{n+1}\}$.  In addition, we
define $\xi=\{\bs{\xi},\xi_{n+1}\}$ such that $\xi$ is isotropic with
respect to the Euclidean scalar product $(\;{,}\,)$ in~$\oC^{n+1}$:
\begin{equation*}
  (\xi,\xi)=\bs{\xi}\ldot\bs{\xi} + \xi_{n+1}^2=0
\end{equation*}
(the other scalar products remain the same as among the
$n$-dimensional vectors above).  With a slight abuse of notation, we
then let $\SPsi$ and $E_i$ denote the screening operators
\begin{equation}\label{n[0]-scr}
  E_i=\oint e^{(a_i,\varphi)},~i=1,\dots,n\,{-}\,1,
  \qquad
  \SPsi=\oint e^{(\psi,\varphi)}.
\end{equation}

\textit{The $n[0]$ realization of $\WW{n}$, denoted by $\WW{n[0]}$, is
  the centralizer of $(E_i)_{i=1,\dots,n-1}$, and~$\SPsi$
  in~$\algV_{\xi} = \{\polP(\dd\varphi)\,e^{m\,(\xi,\varphi)}\mid
  m\,{\in}\,\oZ\}$}.

In what follows, the differential polynomials in $\dd\varphi$ are
expressed through the currents
\begin{gather*}
  A_i=(a_i,\dd\varphi),\qquad Q=(\psi,\dd\varphi),
  \qquad Y=(\xi,\dd\varphi),
\end{gather*}
which have the nonzero operator products (read off from the
matrix~\eqref{Cartan-asym-dressed})
\begin{gather*}
  A_i(z)\,A_{i+1}(w)=\ffrac{-k\,{-}\,n}{(z\,{-}\,w)^2},\quad
  A_i(z)\,A_{i}(w)=\ffrac{2(k\,{+}\,n)}{(z\,{-}\,w)^2},\\
  A_1(z)\,Q(w)=\ffrac{-k\,{-}\,n}{(z\,{-}\,w)^2},\quad
  Q(z)\,Q(w)=\ffrac{1}{(z\,{-}\,w)^2},\quad
  Q(z)\,Y(w)=\ffrac{1}{(z\,{-}\,w)^2}.
\end{gather*}

\medskip

\begin{rem}
  It follows from~\eqref{xi-cond-0} that the $n$-dimensional part
  $\bs{\xi}$ of $\xi$ is given by
  \begin{equation}\label{xi-explicit}
    \bs{\xi}=-\ffrac{1}{n\ell_n(k)}
    \bigl(\bs{a}_{n-1}+2\bs{a}_{n-2}+\dots+(n\,{-}\,1)\bs{a}_{1}+
    n\bs{\psi}\bigr)
  \end{equation}
  and
  \begin{equation*}
    \bs{\xi}\ldot\bs{\xi}=-\ffrac{1}{\ell_n(k)}.
  \end{equation*}
  We also note that the determinant of the
  $(n\,{+}\,1)\,{\times}\,(n\,{+}\,1)$ Gram matrix of the vectors
  $a_{n-1}$, \dots, $a_1$, $\psi$, $\xi$,
  \begin{equation}\label{Gamma}  
    \Gamma_n(k)=
    \mbox{\addtolength{\arraycolsep}{-2pt}
      \small$\displaystyle\begin{pmatrix}
        2(k\,{+}\,n)
        &-k\,{-}\,n&0&\hdotsfor{3}&0&0\\
        -k\,{-}\,n&2(k\,{+}\,n)&-k\,{-}\,n&0&\hdotsfor{2}&0&0\\
        \hdotsfor{8}\\
        \hdotsfor{8}\\
        0&\hdotsfor{2}&0&-k\,{-}\,n&2(k\,{+}\,n)&-k\,{-}\,n&0\\
        0&\hdotsfor{3}&0&-k\,{-}\,n&1&1\\
        0&\hdotsfor{4}&0&1&0
      \end{pmatrix}\!,$}
  \end{equation}
  is given by $-n(k\,{+}\,n)^{n-1}$, and hence these $n\,{+}\,1$
  vectors form a basis in~$\oC^{n+1}$ and are determined uniquely
  modulo a common rotation for $k\,{\in}\,\oC\setminus\{-n\}$. \ 
% (for $k=-\frac{n(n-2)}{n-1}$, the vectors $\psi$, $a_1$, \dots,
% $a_{n-1}$ cannot be chosen in an $n$-dimensional subspace).
  For future use, we note that
  \begin{multline*}
    \Gamma_n(k)^{-1}=\ffrac{1}{n}\times{}\\*
    \addtolength{\arraycolsep}{-2pt}
    \begin{pmatrix}
      \frac{(n-1)\cdot1}{k+n}&\frac{(n-2)\cdot1}{k+n}
      &\frac{(n-3)\cdot1}{k+n}
      &\hdotsfor{4}&\frac{1\cdot1}{k+n}&0&1\\[3pt]
      \frac{(n-2)\cdot1}{k+n}&\frac{(n-2)\cdot2}{k+n}
      &\frac{(n-3)\cdot2}{k+n} &\frac{(n-4)\cdot2}{k+n}&\hdotsfor{3}&
      \frac{1\cdot2}{k+n}&0&2\\[3pt]
      \frac{(n-3)\cdot1}{k+n}&\frac{(n-3)\cdot2}{k+n}
      &\frac{(n-3)\cdot3}{k+n}&\frac{(n-4)\cdot3}{k+n}
      &\frac{(n-5)\cdot3}{k+n}
      &\hdotsfor{2}&\frac{1\cdot3}{k+n}&0&3\\[3pt]
      \frac{(n-4)\cdot1}{k+n}&\frac{(n-4)\cdot2}{k+n}
      &\frac{(n-4)\cdot3}{k+n}&\frac{(n-4)\cdot4}{k+n}
      &\frac{(n-5)\cdot4}{k+n}&\frac{(n-6)\cdot4}{k+n}
      &\dots&\frac{1\cdot4}{k+n}&0&4\\[3pt]
      \hdotsfor{10}\\[3pt]
      \hdotsfor{10}\\[3pt]
      \frac{1\cdot1}{k+n}&\frac{1\cdot2}{k+n}
      &\frac{1\cdot3}{k+n}&\frac{1\cdot4}{k+n}&\hdotsfor{2}
      &\dots&\frac{1\cdot(n-1)}{k+n}&0&n{-}1\\[3pt]
      0&0&\hdotsfor{5}&0&0&n\\[3pt]
      1&2&3&4&\hdotsfor{3}&n{-}1&n&n\ell_n(k)
    \end{pmatrix}.
  \end{multline*}
\end{rem}

\subsubsection{$\boldsymbol{n{[}1{]}}$}\label{sec:n[1]} The
realization $\WW{n[1]}$ corresponds to the Dynkin diagram
\begin{equation*}
  \xymatrix@1{%
    {\bullet}\ar@{-}[r]
%     &{\bullet}\ar@{-}[r]
    &{\bullet}\ar@{-}[];[]+<16pt,0pt>
    &{\dots}
    &{\bullet}\ar@{-}[];[]-<16pt,0pt>
    \ar@{-}[r]
    &{\circ}\ar@{-}[r]
    &{\circ}
  }
\end{equation*}
with two odd roots.  Similarly to the $n[0]$ case, we introduce
vectors $\bs{a}_{n-2}$, \dots, $\bs{a}_1$, $\bs{\psi}_+$,
$\bs{\psi}_-$ in~$\oC^n$ whose Gram matrix is given by the Cartan
matrix ``dressed'' into
\begin{equation*}
  \addtolength{\arraycolsep}{-2pt}
  \mbox{\small$\displaystyle
    \begin{pmatrix}
      2(k\,{+}\,n)\rlap{\kern-90pt$\smash{\begin{array}[t]{r}
            \bs{a}_{n-2}\\\bs{a}_{n-3}\\ \\ \\ \bs{a}_1\\
            \bs{\psi}_+\\\bs{\psi}_-
          \end{array}\rule[-94pt]{.5pt}{100pt}}$}
      &-k\,{-}\,n  &0   &\hdotsfor{3}&0 %&0
      \\
      -k\,{-}\,n          &2(k\,{+}\,n)&-k\,{-}\,n&0&\hdotsfor{2}&0 %&0
      \\
      \hdotsfor{7}\\
      \hdotsfor{7}\\
      0&\hdotsfor{1}&0     &-k\,{-}\,n&2(k\,{+}\,n)&-k\,{-}\,n &0 %&0
      \\
      0&\hdotsfor{2}&0     &-k\,{-}\,n&1     &k\,{+}\,n\,{-}\,1 %&1
      \\
      0&\hdotsfor{3}&0&    k\,{+}\,n\,{-}\,1&1     %&-1
%     \\
%     0&\hdotsfor{3}&0     &1&-1&0
    \end{pmatrix}.$}\kern-40pt
\end{equation*}
We also introduce an $n$-tuple of scalar fields~$\bs{\varphi}$ and
define the screenings
\begin{equation*}
  E_i=\oint e^{\bs{a}_i\ldot\bs{\varphi}},~i=1,\dots,n-2,
  \qquad
  \SPsi_+=\oint e^{\bs{\psi}_+\ldot\bs{\varphi}},\quad
  \SPsi_-=\oint e^{\bs{\psi}_-\ldot\bs{\varphi}},
\end{equation*}
where the dot denotes the Euclidean scalar product in~$\oC^n$.  These
screenings represent $\nqSSL{n}{1}$.  We next define the vector
$\bs{\xi}$ by its products with $\bs{a}_{n-2}$, \dots, $\bs{a}_1$,
$\bs{\psi}_+$, $\bs{\psi}_-$,
\begin{equation}\label{xi-cond-1}
  \begin{split}    
    \bs{\xi}\ldot\bs{a}_i &=0,\quad i=1,\dots,n-2,\\
    \bs{\xi}\ldot\bs{\psi}_+&=1,\\
    \bs{\xi}\ldot\bs{\psi}_-&=-1.
  \end{split}
\end{equation}
The $n[1]$ realization of $\bWW{n}$, denoted by $\bWW{n[1]}$, is the
centralizer of $(E_i)_{i=1,\dots,n-2}$, $\SPsi_+$, and $\SPsi_-$
in~$\bar{\algV}_{\bs{\xi}} =
\{\polP(\dd\bs{\varphi})\,e^{m\,\bs{\xi}\ldot\bs{\varphi}}\mid
m\,{\in}\,\oZ\}$.

To construct $\WW{n[1]}$, we embed $\oC^n$ in $\oC^{n+1}$ as the
coordinate hyperplane, let $a_{n-2}$, \dots, $a_1$ $\psi_+$, $\psi_-$
denote the respective images of $\bs{a}_{n-2}$, \dots, $\bs{a}_1$,
$\bs{\psi}_+$, $\bs{\psi}_-$, and extend the free scalar fields
as~$\dd\varphi=\{\dd\bs{\varphi},\dd\varphi_{n+1}\}$.  In addition, we
define $\xi=\{\bs{\xi},\xi_{n+1}\}$ such that $\xi$ is isotropic with
respect to the Euclidean scalar product $(\;{,}\,)$ in~$\oC^{n+1}$
(the other scalar products remain the same as among the
$n$-dimensional vectors above).  We again use $\SPsi_-$, $\SPsi_+$,
and $E_i$ to denote the screenings
\begin{equation*}
  E_i=\oint e^{(a_i,\varphi)},~i=1,\dots,n-2,\quad
  \SPsi_+=\oint e^{(\psi_+,\varphi)},\quad
  \SPsi_-=\oint e^{(\psi_-,\varphi)}.
\end{equation*}
The screenings $\SPsi_+$ and $\SPsi_-$ are said to be
\textit{fermionic} and the other screenings \textit{bosonic}.

\textit{The $n[1]$ realization $\WW{n[1]}$ of~$\WW{n}$ is the
  centralizer of $(E_i)_{i=1,\dots,n-2}$, $\SPsi_+$, $\SPsi_-$ in
  $\algV_{\xi} = \{\polP(\dd\varphi)\,e^{m\,(\xi,\varphi)}\mid
  m\,{\in}\,\oZ\}$}.  As in the maximally asymmetric case, we express
differential polynomials in $\dd\varphi$ through the currents
\begin{gather*}
  A_i=(a_i,\dd\varphi),\qquad Q_{\pm}=(\psi_{\pm},\dd\varphi),
  \qquad Y=(\xi,\dd\varphi),
\end{gather*}
whose OPEs are determined by the matrix above and scalar
products~\eqref{xi-cond-1},
\begin{gather*}
  A_i(z)\,A_{i+1}(w)=\ffrac{-k\,{-}\,n}{(z\,{-}\,w)^2},\quad
  A_i(z)\,A_{i}(w)=\ffrac{2(k\,{+}\,n)}{(z\,{-}\,w)^2},\\
  A_1(z)\,Q_+(w)=\ffrac{-k\,{-}\,n}{(z\,{-}\,w)^2},\quad
  Q_+(z)\,Q_-(w)=\ffrac{k\,{+}\,n\,{-}\,1}{(z\,{-}\,w)^2},\\
  Q_{\pm}(z)\,Q_{\pm}(w)=\ffrac{1}{(z\,{-}\,w)^2},\quad
  Q_{\pm}(z)\,Y(w)=\ffrac{\pm1}{(z\,{-}\,w)^2}.
\end{gather*}

\medskip

\begin{rem}
  It follows that
  \begin{equation*}
    \bs{\xi}=-\ffrac{1}{n\ell_n(k)}\bigl(
    \bs{a}_{n-2}\,{+}\,2\bs{a}_{n-3}\,{+}\,\dots\,{+}\,(n-2)\bs{a}_1
   \,{+}\,(n\,{-}\,1)\bs{\psi}_+
   \,{-}\,\bs{\psi}_-\bigr)
  \end{equation*}
  and
  \begin{equation*}
    \bs{\xi}\ldot\bs{\xi}=-\ffrac{1}{\ell_n(k)}.
  \end{equation*}
  The determinant of the $(n\,{+}\,1)\times(n\,{+}\,1)$ Gram matrix of
  the vectors $a_{n-2}$, \dots, $a_1$, $\psi_+$, $\psi_-$, $\xi$ is
  given by $-n(k\,{+}\,n)^{n-1}$, and these $n\,{+}\,1$ vectors
  therefore form a basis in~$\oC^{n+1}$ and are determined uniquely
  modulo a common rotation for all $k\,{\in}\,\oC\setminus\{-n\}$.
\end{rem}

\subsubsection{$\boldsymbol{n{[}m{]}}$}\label{sec:n[m]}
The subsequent $n[m]$ realizations correspond to the Dynkin diagrams
\begin{equation*}
  \xymatrix@1{%
    {\bullet}\ar@{-}[r]
%     &{\bullet}\ar@{-}[r]
    &{\bullet}\ar@{-}[];[]+<16pt,0pt>
    &{\dots}
    &{\bullet}\ar@{-}[];[]-<16pt,0pt>
    \ar@{-}[r]
    &{\circ}\ar@{-}[r]
    &{\circ}\ar@{-}[r]
    &{\bullet}\ar@{-}[];[]+<16pt,0pt>
    &{\dots}
    &{\bullet}\ar@{-}[];[]-<16pt,0pt>
  },
\end{equation*}
with $m\,{-}\,1$ even roots to the right of the odd roots.  We can
restrict $m$ to $0\,{\leq}\, m\,{\leq}\,[\frac{n}{2}]$, because taking
$m>[\frac{n}{2}]$ amounts to applying an automorphism (inducing
$\cE\leftrightarrow\cF$ on the algebra) to the realization where $m$
is replaced by $n\,{-}\,m$.

We do not repeat the definition of $\bWW{n[m]}$ in terms of $n$ scalar
fields and proceed to defining $\WW{n[m]}$ in terms of $n\,{+}\,1$
scalar fields.  For this, we introduce $(n\,{+}\,1)$-dimensional
vectors $a_{n-m-1}$, \dots, $a_{1}$, $\psi_+$, $\psi_-$, $a_{-1}$,
\dots, $a_{-m\,{+}\,1}$, $\xi$ whose Gram matrix (with the determinant
$-n(k\,{+}\,n)^{n-1}$) is given by
\begin{equation*}  
  \mbox{\footnotesize$\displaystyle
    \setcounter{MaxMatrixCols}{15}
    \addtolength{\arraycolsep}{-2pt}
    \kern60pt\begin{pmatrix}
      \uline{2K}\rlap{\kern-80pt$\smash{\begin{array}[t]{r}
            a_{n-m-1}\\a_{n-m-2}\\ \smash{\vdots} \\a_1\\\psi_+\\
            \psi_-\\a_{-1}\\ \smash{\vdots}\\a_{-2}\\a_{-m+1}\\\xi
          \end{array}\rule[-140pt]{.5pt}{146pt}}$}
      &-K  &0   &\hdotsfor{8}&0&0 \\
      -K        &\uline{2K}&-K&0&\hdotsfor{7}&0&0 \\
      \hdotsfor{13}\\
      0&\hdotsfor{1}&0     &-K&\uline{2K}&-K &0&\hdotsfor{4}&0&0 \\
      0&\hdotsfor{2}&0     &-K&\uline{1}&K-1&0&\hdotsfor{3}&0&1 \\
      0&\hdotsfor{3}&0&    K-1&\uline{1}&-K&0&\hdotsfor{2}&0&-1 \\
      0&\hdotsfor{4}&0&    -K& \uline{2K}&-K&0&\hdotsfor{1}&0&0 \\
      \hdotsfor{13}\\
      0&\hdotsfor{7}&0&    -K& \uline{2K}&-K&0 \\
      0&\hdotsfor{8}&0&    -K& \uline{2K}&0 \\
      0&\hdotsfor{3}&0     &1&-1&0&\hdotsfor{2}&0&0&\uline{0}
    \end{pmatrix},$}
\end{equation*}
where $K=k\,{+}\,n$, the leftmost column indicates labeling of rows, and
diagonal elements are underlined to guide the eye.  
%% The reader can easily truncate this to $n$-dimensional vectors and
%% thus reconstruct the definition of $\bWW{n[m]}$.
%% We note that the $n$-dimensional part $\bs{\xi}$ of $\xi$ is
%% determined by the scalar products read off from the above matrix
%% as
%% \begin{multline*}
%%   \bs{\xi}=  -\ffrac{1}{n\ell_n(k)}
%%   \bigl(\bs{a}_{n-m-1}
%%   +2\bs{a}_{n-m-2} + \dots + (n-m-1)\bs{a}_1
%%   +(n-m)\bs{\psi}_+\\*
%%   -\bs{a}_{-1}-2\bs{a}_{-2}
%%   -\dots-(m-1)\bs{a}_{-m+1} - m\bs{\psi}_-\bigr).
%% \end{multline*}
The screenings that represent $\nqSSL{n}{1}$ are given by
\begin{multline*}
  E_i=\oint e^{(a_i,\varphi)},~i=1,\dots,n\,{-}\,m\,{-}\,1,\\
  \SPsi_+=\oint e^{(\psi_+,\varphi)},\quad
  \SPsi_-=\oint e^{(\psi_-,\varphi)},\\
  E_i=\oint e^{(a_i,\varphi)},~i=-1,\dots,-m+1.
\end{multline*}
\textit{The $n[m]$ realization of $\WW{n}$, denoted by $\WW{n[m]}$, is
  the centralizer of these operators in $\algV_{\xi} =
  \{\polP(\dd\varphi)\,e^{m\,(\xi,\varphi)}\mid m\,{\in}\,\oZ\}$}.  As
before, $E_i$ are said to be bosonic and $\SPsi_{\pm}$ fermionic
screenings.

The information contained in the Gram matrix above can be conveniently
reexpressed as a ``rigged'' Dynkin diagram for $\SSL{n}{1}$,
\begin{equation}\label{rigged}
  \rule[-24pt]{0pt}{36pt}
  \xymatrix{%
    *{\;\bullet\;}\ar@{-}[0,1]^{-K}
    \ar@(ld,rd)@{-}[]|{2K}
    &*{\;\bullet\;}\ar@{-}[0,1]^{-K}
    \ar@(ld,rd)@{-}[]|{2K}
    &{\dots}\ar@{-}[0,1]^{-K}
    &*{\;\bullet\;}\ar@{-}[0,1]^{-K}
    \ar@(ld,rd)@{-}[]|{2K}
    &*{\;\circ\;}\ar@{-}[0,1]^{K-1}
    \ar@(ld,rd)@{-}[]|{1}
    \ar@{}[]+<0pt,20pt>|(.4){\boldsymbol{+}}
    &*{\;\circ\;}\ar@{-}[0,1]^{-K}
    \ar@(ld,rd)@{-}[]|{1}
    \ar@{}[]+<0pt,20pt>|(.4){\boldsymbol{-}}
    &*{\;\bullet\;}\ar@{-}[0,1]^{-K}
    \ar@(ld,rd)@{-}[]|{2K}
    &{\dots}\ar@{-}[0,1]^{-K}
    &*{\;\bullet\;}
    \ar@(ld,rd)@{-}[]|{2K}\\
    *{\!\mystrut{A_{n-m-1}}\!\!\!\!}
    &*{\!\!\!\!\mystrut{A_{n-m-2}}\!\!\!\!}
    &*{\!\!\mystrut{\dots}\!\!}
    &*{\mystrut{A_{1}}}
    &*{\mystrut{Q_+}}
    &*{\mystrut{Q_-}}
    &*{\mystrut{A_{-1}}}
    &*{\!\!\mystrut{\dots}\!\!\!\!}
    &*{\!\!\!\!\mystrut{A_{1-m}}\!\!}
    }
\end{equation}

\medskip

\noindent
where $K=k\,{+}\,n$.  The vertices of the diagram are assigned the
dimension-$1$ currents $A_{n-m-1}(z)$, \dots, $A_{1}(z)$, $Q_+(z)$,
$Q_-(z)$, $A_{-1}(z)$, \dots, $A_{-m+1}(z)$ as indicated.
%% $A_{n-m-1}$, \dots, $A_{1}$ (the left part of the diagram, from left
%% to right) and $A_{-1}$, \dots, $A_{-m+1}$ (the right part of the
%% diagram, from left to right), $A_i=(a_i,\dd\varphi)$.  The open dots
%% are assigned the currents $Q_+$ and $Q_-$.  
The labels $2K$, $-K$, and $1$ at the links mean that these currents
have the OPEs
\begin{alignat*}{2}
  A_i(z)\,A_i(w)&=\mfrac{2K}{(z\,{-}\,w)^2},&\qquad
  A_i(z)\,A_{i+1}(w)&=\mfrac{-K}{(z\,{-}\,w)^2},\\
  Q_{\pm}(z)\,Q_{\pm}(w)&=\mfrac{1}{(z\,{-}\,w)^2},&
  Q_{+}(z)\,Q_{-}(w)&=\mfrac{K-1}{(z\,{-}\,w)^2},\\
  A_{1}(z)\,Q_{+}(w)&=\mfrac{-K}{(z\,{-}\,w)^2},&
  A_{-1}(z)\,Q_{-}(w)&=\mfrac{-K}{(z\,{-}\,w)^2}.
\end{alignat*}
In addition, for the current $ Y=(\xi,\dd\varphi)$, which is not
associated with a vertex in the diagram, we have the nonzero OPEs (as
indicated by the ${+}$ and ${-}$ signs),
\begin{equation*}
  Q_{\pm}(z)\,Y(w)=\mfrac{\pm1}{(z\,{-}\,w)^2}.
\end{equation*}

\subsection{Centralizer of the screenings: the ``step-back''
  strategy}\label{sec:strategy} We now describe the strategy that
leads to the construction of the $\WW{n[m]}$ generators in the
centralizer of the corresponding screenings.  The contents of this
subsection is not a proof, but because it motivates our construction,
we hope that it can be useful in finding centralizers of some other
quantum supergroups (a necessary condition is the existence of
fermionic screening(s)).  The crucial point is the conditions on the
vector~$\xi$ (see~\eqref{xi-cond-0}, \eqref{xi-cond-1}, and similar
conditions read off from the matrix in~\bref{sec:n[m]}), chosen from
quantum group considerations.

Quite generally, we recall that the action of screening operators
$S_1,\dots,S_N\,{\in}\,\qnilp$ on an operator $X=\polP'V(p',z)$ in a
lattice vertex-operator algebra (where $V(p',z)$ is a vertex with
momentum~$p'$ and $\polP'$ is a differential polynomial) generically
gives nonlocal expressions which represent elements of a module
$\qmodK$ over the quantum group~$\qalg$.  This module can be either a
Verma module or (which is most often the case) some of its quotients.
But whenever a singular vector, e.g., $v=S_1\dots S_r X$ (where we
write a monomial expression for simplicity) occurs in $\qmodK$, the
corresponding field is \textit{local}\,---\,a descendant of the
``shifted'' vertex $V(p,z)$, where the momentum $p$ differs from $p'$
by the sum of the momenta of the relevant screenings,
$p=p'+a_1+\dots+a_r$.  We are interested in the case where a singular
vector in $\qmodK$ generates a \textit{one-dimensional submodule}.
Then the corresponding \textit{local} field is necessarily \textit{in
  the centralizer} of the quantum group.  To construct a local field
of the form $\polP\,V(p,z)$, we must then start with an appropriate
$p'=p-a_1-\dots-a_r$.

That is, we seek a quantum group module~$\qmodK$ with a singular
vector $v$ that generates a \textit{one-dimensional submodule}
in~$\qmodK$ (the singular vector must therefore have weight~$0$).
Such singular vectors (and hence the corresponding local fields) can
indeed be found for the $\SSL{n}{1}$ root system (the existence of odd
roots is a necessary condition).  We now use this method to construct
fields in the centralizer of $\qSSL{n}{1}$.

\subsubsection{Maximally asymmetric realization}
We first consider the maximally asymmetric realization, corresponding
to the $\SSL{n}{1}$ simple root system described in~\bref{sec:n[0]}.
We must find a Verma module quotient with a singular vector generating
a one-dimensional submodule.  For this, we take the
$\qSSL{n}{1}$-module induced from the one-dimensional representation
of the parabolic subalgebra generated by $\Univ_q\GL{n}$ and the
fermionic simple root generator.  The module is then isomorphic to
$\wwedge{}\kern-2pt\raisebox{6pt}{$\scriptstyle\bullet$}\kern2pt\oC^n$
as a vector space.  Its weight diagram is shown in Fig.~1 for $n=3$
(the upper-left part, with the highest-weight vector at the top).
\begin{SCfigure}[1][tb]
  \raisebox{11\baselineskip}{
  \xymatrix@=16pt{
    &{\circ}\\
    {\circ}&{\circ}&{\circ}\\
    {\circ}&{\circ}&{\circ}\\
    &{\bullet}&&&&{\circ}\ar[]+<-12pt,0pt>;[0,-4]+<12pt,0pt>\\
    &&&&{\bullet}&{\bullet}&{\bullet}\\
    &&&&{\bullet}&{\bullet}&{\bullet}\\
    &&&&&{\bullet}
  }}
  \caption[Mapping between two $\SSL31$ modules]{\small Mapping
    between two $\SSL31$ modules whose highest-weight vectors are
    annih\-il\-ated by the~$\SL3$ subalgebra.  Filled dots denote
    elements of \textit{submodules}.  Those in the right module
    ($M_3(0)$) vanish under the mapping into the left module
    ($M_3(-2)$).}\label{fig:sl31q}
\end{SCfigure}
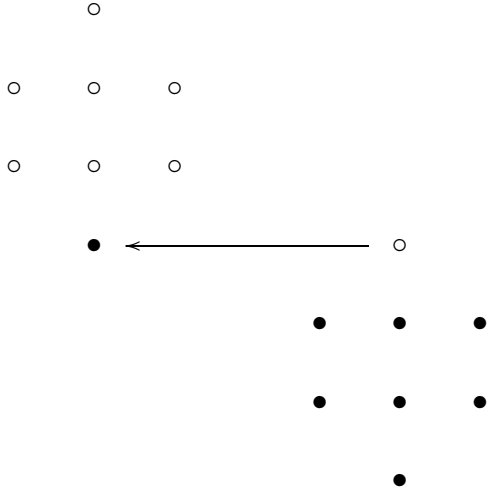
For generic $q$, such a $\qSSL{n}{1}$-module is a deformation of an
$\SSL{n}{1}$-module and is isomorphic to the latter as a vector space,
and we proceed with the corresponding $\SSL{n}{1}$-module for
simplicity.  Such $\SSL{n}{1}$-modules depend on a single parameter,
the eigenvalue $\alpha$ of the $\GL{n}$ Cartan generator that is not
in~$\SL{n}$.  We normalize this generator as
\begin{equation*}
  h_0=
  \begin{pmatrix}
    \frac{1}{n}\\
    &\ddots\\
    &&\frac{1}{n}\\
    &&&1
  \end{pmatrix}
\end{equation*}
and let $\ket{\alpha}$ be the highest-weight vector of the above
module with $h_0\ket{\alpha} = \alpha\ket{\alpha}$; the module is
denoted by~$M_n(\alpha)$.  \textit{The vector at the bottom of the
  weight diagram of~$M_n(\alpha)$ is then singular if and only
  if~$\alpha=-(n\,{-}\,1)$}.

In particular, the vector at the bottom of the weight diagram
of~$M_n(1\,{-}\,n)$ is annihilated by $h_0$, and the highest-weight
vector of the $M_{n}(0)$ module is then mapped onto this singular
vector, see~Fig.~1.  This picture is preserved under deformation to
$\qSSL{n}{1}$.  By the correspondence between modules over the
quantum-group and the vertex operator algebra, the singular vector in
(the deformation of)~$M_n(1\,{-}\,n)$ determines an intertwining
operator between the corresponding $\WW{n[0]}$ representations,
realized in a sum of Fock spaces; this intertwining operator is
morally the product of the $n$ fermionic screenings $\SPsi$, $\SPsi'$,
$ \dots$, $\SPsi^{(n-1)}$ obtained by the action of the bosonic
screenings on~$\SPsi$,
\begin{gather*}
  \SPsi\SPsi'\dots\SPsi^{(n-1)} : F_{\bs{p}'}\to F_{\bs{p}},
\end{gather*}
where $F_{\bs{p}}$ is the module generated from
$e^{\bs{p}\ldot\bs{\varphi}}$.  The image of the highest-weight vector
in~$F_{\bs{p}'}$ is a descendant of the highest-weight vector
in~$F_{\bs{p}}$,
\begin{gather}\label{F-to-F}
  (\SPsi\SPsi'\dots\SPsi^{(n-1)})(e^{\bs{p}'\ldot\bs{\varphi}})
  =\polP_n\,e^{\bs{p}\ldot\bs{\varphi}},
\end{gather}
where $\polP_n$ is a degree-$n$ differential polynomial in
$\dd\bs{\varphi}= (\dd\varphi_1,\dots,\dd\varphi_n)$.  The difference
$\bs{p}-\bs{p}'$ must be equal to the sum of the momenta of $\SPsi$,
$\SPsi'$, \dots, $\SPsi^{(n-1)}$, which are $\bs{\psi}$ for the
fermionic screening in~\eqref{n[0]-scr} and
$\bs{\psi}+\bs{a}_1+\dots+\bs{a}_{i}$ for each of the other fermionic
screenings obtained by dressing $\SPsi$ with the bosonic screenings.
Therefore,
\begin{equation}\label{step-back}
  \bs{p} = \bs{p}' + n\bs{\psi} + (n\,{-}\,1)\bs{a}_1
  + (n-2)\bs{a}_2 + \dots + \bs{a}_{n-1}
  =\bs{p}'-n\ell_n(k)\bs{\xi},
\end{equation}
with $\bs{\xi}$ defined in~\eqref{xi-explicit}.

Next, to ensure that the module is induced from the one-dimensional
representation of the parabolic subalgebra, we require that
\begin{equation*}
  \bs{p}\ldot\bs{a}_i=0,\quad i=1,\dots,n\,{-}\,1
\end{equation*}
(with the scalar products among $\bs{a}_i$ and $\bs{\psi}$ given
in~\bref{sec:n[0]}, this is equivalent to $\bs{p}'\ldot\bs{a}_i=0$).
Therefore, $[E_i,e^{\bs{p}\ldot\bs{\varphi}}]=0$ and
$[E_i,e^{\bs{p'}\ldot\bs{\varphi}}]=0$.  Further, the vector
in~$\qmodK$ represented by the bottom dot in the upper-left module in
Fig.~1 is singular, and the module is therefore (the quantum
deformation of) $M_n(1\,{-}\,n)$, if
%% we choose $\bs{p}$ such that
\begin{equation*}
  \bs{p}\ldot\bs{\psi}=-1.
\end{equation*}

With these conditions satisfied, applying
$(\SPsi\SPsi'\dots\SPsi^{(n-1)})$ to $e^{\bs{p}'\ldot\bs{\varphi}(z)}$
gives a level-$n$ descendant $\polP_n\,e^{\bs{p}\ldot\bs{\varphi}(z)}$
of $e^{\bs{p}\ldot\bs{\varphi}(z)}$ that necessarily commutes with the
screenings.\footnote{More generally, we expect that for generic $k$,
  the entire centralizer of $\qSSL{n}{1}$ in the sector with momentum
  $\bs{p}$ is given by $(\SPsi\SPsi'\dots\SPsi^{(n-1)})\,(\polP'\,
  e^{\bs{p}'\ldot\bs{\varphi}(z)})$, where $\polP'$ is a differential
  polynomial such that $[E_i,\polP' e^{\bs{p'}\ldot\bs{\varphi}}]=0$.}
It follows that $\bs{p}=-\bs{\xi}$, see~\eqref{xi-explicit}.  The
currents
\begin{align*}
  \Bar{\cE}_{n[0]}(z)={}&e^{\bs{\xi}\ldot\bs{\varphi}(z)},\\
  \Bar{\cF}_{n[0]}(z)={}&{-\polP_n}\,e^{-\bs{\xi}\ldot\bs{\varphi}(z)}
\end{align*}
then generate $\bWW{n[0]}$.  To recover $\WW{n[0]}$, it remains to
introduce an additional free scalar field and embed the above
$n$-dimensional vectors in $\oC^{n+1}$ as $a_i=\{\bs{a}_i,0\}$,
$\psi=\{\bs{\psi},0\}$, and $\xi=\{\bs{\xi},\xi_{n+1}\}$ with
\textit{isotropic}~$\xi$, as in~\bref{sec:n[0]}.  We introduce the
currents
\begin{equation}\label{A-currents}
  A_i=(a_i,\dd\varphi),\quad Q=(\psi,\dd\varphi)
\end{equation}
and the scalar field
\begin{equation*}
  \YXi=(\xi,\varphi)
\end{equation*}
such that
\begin{gather}\label{Y-d-Xi}
  Y=\dd\,\YXi=(\xi,\dd\varphi).
\end{gather}
{}From the scalar product $(\psi,\xi)=1$, we then have the OPEs
\begin{gather*}%%\label{Q-Xi}
  Q(z)\,e^{\pm\YXi(w)}
  =\mfrac{\pm1}{z\,{-}\,w}\,e^{\pm\YXi(w)}.
\end{gather*}
In the ``maximally asymmetric'' realization, the $\WW{n[0]}$-currents
have the form
\begin{gather*}%%%\label{X-n[0]}
  \begin{aligned}
    \cE_{n[0]}(z) ={}& e^{\YXi}(z),\\
    \cF_{n[0]}(z) ={}& {-}\polP_{n}(A_{n-1},\dots,A_1,Q)
    \,e^{-\YXi(z)},
  \end{aligned}
\end{gather*}
where $\polP_{n}(A_{n-1},\dots,A_1,Q)$ is the polynomial
in~\eqref{F-to-F} expressed through $A_{n-1}$, \dots, $A_1$,~$Q$.

\subsubsection{Other realizations}
In other realizations, both the $\cE$ and $\cF$ currents are given by
a normal-ordered product of an exponential and a differential
polynomial.  We recall that in the well-known ``symmetric''
bosonization of $\hSL2$ (see~\bref{sec:sl2}), the $\cE$ and $\cF$
currents follow by the action of the corresponding fermionic screening
and are therefore given by order-$1$ polynomials in front of the
exponentials, see~\eqref{sl2-sym}.  In the general case, the simplest
singular vectors given by the action of a single $\qSSL{2}{1}$
generator are replaced with singular vectors in the appropriate
$\qSSL{n}{1}$ modules.

A ``step-back'' strategy similar to the one used in the maximally
asymmetric case also involves $2^n$-dimensional $\qSSL{n}{1}$-modules;
the classical $\SSL{n}{1}$-analogue of every such module can be viewed
as an $M_n(1\,{-}\,n)$ module turned on its side.  This is illustrated
in Fig.~\figsliv\ for $n=4$.
\begin{figure}%[20][tb]
  \mbox{\xymatrix@=16pt{
    &&{}\ar@{}|{\mbox{\normalsize$\circ$}}[r]&\\
    &{\circ}&{\circ}&{\circ}&{\diamond\smash{%
        \makebox[0pt]{\raisebox{10pt}{$s$}}}}\\
    {\circ}&{\circ}&{\circ}&{\diamond}&{\diamond}&{\ddiamond}\\
    &{\circ\smash{%
        \makebox[0pt]{\kern-14pt\raisebox{-12pt}{$s'$}}}}&
    {\diamond}&{\ddiamond}&{\ddiamond}\\    
    &&{}\ar@{}|{\mbox{\normalsize$\bullet$}}[r]&&&&}
  \quad%%\qquad\qquad\qquad\qquad\qquad
      \raisebox{-6.5\baselineskip}{%
        \xymatrix@=4pt{\bullet&\bullet&\bullet&\circ
      &&
      \bullet&\bullet&\circ&\circ
      \\
      {}&{}\ar@{}[r]|{(a)}&{}&{}
      &&
      {}&{}\ar@{}[r]|{(b)}&{}&{}}}}
  \caption[A ``matrioshka'' arrangement of modules]{\small A
    ``matrioshka'' arrangement of special $\SSL{n}{1}$-modules
    (\textit{left}).  Considered as generated from the top vector by
    the simple root generators corresponding to the Dynkin
    diagram~(\textit{a}), the $2^4$ states represent the
    $\protect\SSL41$-module $M_4(-3)$.  The diamonds show the
    $\protect\SSL31$-module $M_3(-2)$, as in the left part of
    Fig.~\figsliii.  The filled dot is a singular vector in both
    $M_4(-3)$ and $M_3(-2)$.  The bigger diamonds (together with
    $\bullet$, which is again a singular vector) form the weight
    diagram of the $\protect\SSL21$-module $M_2(-1)$.  With the simple
    root generators corresponding to the Dynkin diagram~(\textit{b}),
    the $2^4$-dimensional $\SSL41$ module is generated from the
    state~$s$.}\label{fig:sl41q}
\end{figure}
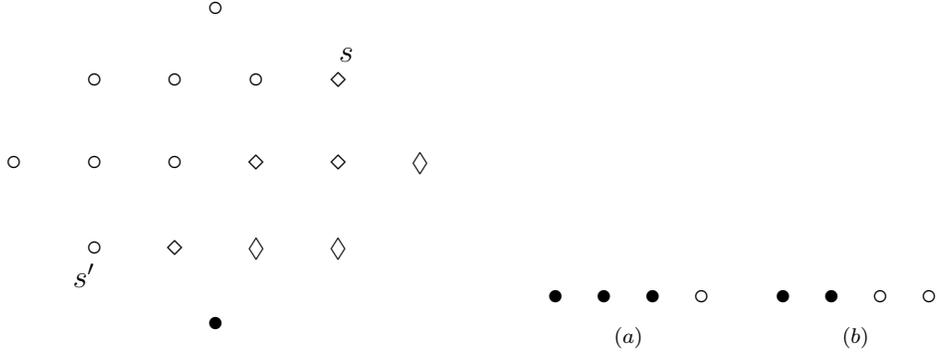
If the module is viewed as generated from the top vector, it is
$M_4(-3)$.  In this case, the simple root generators correspond to the
Dynkin diagram~(\textit{a}).  But with the simple root generators
corresponding to the Dynkin diagram~(\textit{b}), the
$2^4$-dimensional $\SSL41$ module is generated from the vector~$s$
east-south-east of the top one.  The vector at the bottom ``in the
$M_4(-3)$ coordinates'' is then again singular and generastes a
one-dimensional submodule.  Moreover, restricting to the $\SSL31$
subalgebra (corresponding to the first three nodes in
diagram~(\textit{b})), we see that the $\SSL31$-module generated
from~$s$ (shown with diamonds in Fig.~\figsliv) is $M_3(-2)$, as in
the left-hand part of Fig.~\figsliii.  The filled dot is also a
singular vector in $M_3(-2)$, with a one-dimensional submodule (it is
also clear how the $\SSL21$-module $M_2(-1)$ fits this picture (bigger
diamonds).  Construction of the singular vector from the~$s$ state,
actually \textit{via the $\SSL{3}{1}$ generators}, then translates
into a formula for the $\cF$ current in the centralizer of the
screenings.

Moreover, the $\cE$ current in the centralizer of the screenings then
follows by dressing the corresponding exponential with a single simple
root generator; in Fig.~\figsliv, the corresponding module should be
viewed as generated from the state~$s'$, at ``distance~$1$'' from the
singular vector with a $1$-dimensional submodule.

This hierarchy of the special $\SSL{n}{1}$-modules under consideration
(in fact, of their $\qSSL{n}{1}$-deformations) shows that singular
vectors generating one-dimensional submodules can be constructed via
the same step-back mechanism as in the maximally asymmetric case
described above, but applied to a subset of simple root generators.
The singular vector is in a sense ``the same'' for the different root
systems in the $\SSL{n}{1}$ algebras with different~$n$, only
constructed from differently chosen highest-weight vectors by the
appropriate set of simple root generators.  Accordingly, the $\cE$ and
$\cF$ currents in the centralizer of the screenings are then obtained
by dressing the corresponding exponential with \textit{complementary}
sets of the simple root generators, $E_i$,
$i=1,\dots,n\,{-}\,m\,{-}\,1$, and $\SPsi_+$ for $\cF_{n[m]}$, and
$E_i$, $i=-1,\dots,-m+1$ and $\SPsi_-$ for $\cE_{n[m]}$, and
therefore,
\begin{gather*}%%%\label{X-n[0]}
  \begin{aligned}
    \cE_{n[m]}(z) ={}& \polP^\dagger_{m}(A_{-1},\dots,A_{-m+1},Q_-)
    e^{\YXi}(z),\\
    \cF_{n[m]}(z) ={}& (-)^{m+1}\polP_{n-m}(A_{n-m-1},\dots,A_1,Q_+)
    \,e^{-\YXi(z)}
  \end{aligned}
\end{gather*}
with differential polynomials $\polP^\dagger_{m}$ and $\polP_{n-m}$ of
the respective order $m$ and $n\,{-}\,m$ (normal ordering in the
right-hand sides is understood).

\subsection{Centralizer of $\qSSL{n}{1}$: recursion
  relations}\label{sec:recursion} Parallel to the ``matrioshka''
arrangement of the $\SSL{n}{1}$ modules described in the previous
subsection, there exist recursion relations between the differential
polynomials entering the different realizations of the $\WW{n}$
algebras with different~$n$.  In considering these, we must be very
precise about notation.

\subsubsection{Notational chores}\label{sec:chores} We write
$\cX_{n[m]}^{(k)}$ for generators
$\cX_{n[m]}=(\cE_{n[m]},\cH_{n[m]},\cF_{n[m]})$ of $\WW{n[m]}(k)$
whenever we need to indicate the level~$k$ of the algebra.  The
differential polynomials $\polP^\dagger_{n[m]}$ and $\polP_{n[m]}$ in
$A_i$ and $Q_{\pm}$ entering these generators depend on~$k$
explicitly, which we indicate by writing them
as~$\polP^{\dagger(k)}_{n[m]}$ and~$\polP^{(k)}_{n[m]}$.  Moreover,
the free fields introduced in~\bref{sec:n[m]} also bear an implicit
dependence on~$k$ and $n$, in fact on $k\,{+}\,n$ involved in their
OPEs (where we recall that $K=k\,{+}\,n$).  When we need to be very
precise, we write these free fields as $A_i^{[k+n]}$ and
$Q_{\pm}^{[k+n]}$, and similarly use the notation $E_i^{[k+n]}$ and
$\SPsi_\pm^{[k+n]}$ for the screenings.  In the $n[m]$ realization of
$\WW{n}(k)$, we can thus write the generators, most generally, as
(see~\eqref{A-currents}--\eqref{Y-d-Xi} for the definition of the
currents)
\begin{gather}\label{p-a-priori}
  \begin{split}
    \cE^{(k)}_{n[m]}(z)&=
    \polP^{\dagger(k)}_{n[m]}
    (A_\bullet^{[k+n]},Q_+^{[k+n]},Q_-^{[k+n]})\,
    e^{\YXi(z)},\\
    \cF^{(k)}_{n[m]}(z)&= (-1)^{m+1}
    \polP^{(k)}_{n[m]}(A_\bullet^{[k+n]},Q_+^{[k+n]},Q_-^{[k+n]})
    \,e^{-\YXi(z)},
  \end{split}
\end{gather}
where $\polP^{\dagger(k)}_{n[m]}$ and $\polP^{(k)}_{n[m]}$ are
differential polynomials of the respective degrees $m$ and $n-m$ and
$A_\bullet$ stands for the appropriate collection of the $A_i$
currents (see~\bref{sec:n[m]}).  The conventional sign factor is
chosen for future convenience.  By definition,
$\polP^{\dagger}_{n[0]}=\polP^{\phantom{\dagger}}_{n[n]}=1$.  For
$m=0$ (the maximally asymmetric realization), $Q_-$ does not enter and
$Q_+$ is identified with~$Q$ in~\bref{sec:n[0]}.  

\begin{example}\label{ex:n=2}From the well-known three-boson
  realizations of $\hSL2$, see~\bref{sec:sl2}, we have that the lowest
  $\polP$ polynomials are given by
  \begin{equation}\label{rec-bc}
    \begin{split}
      \polP^{\dagger(k)}_{2[1]}(Q)&=\polP^{(k)}_{2[1]}(Q)=Q,\\
      \polP^{(k)}_{2[0]}(A,Q) &= A Q + Q Q + (k + 1) \dd Q.
    \end{split}
  \end{equation}
\end{example}

\subsubsection{Field identifications} It follows that the free fields
entering a realization of the $\WW{n-1}(k\,{+}\,1)$ algebra can be
\textit{identified} with a subset of the fields involved in a
realization of~$\WW{n}(k)$.  We can therefore consider the
``universal'' sets of the currents
\begin{gather*}
  \dots,A^{[\kappa]}_2,A^{[\kappa]}_1,
  Q_+^{[\kappa]}, Q_-^{[\kappa]},
  A^{[\kappa]}_{-1},A^{[\kappa]}_{-2},\dots
\end{gather*}
and the screenings 
%% $(E^{[\kappa]}_i)_{i\,{\in}\,\oZ}$ and $\SPsi_\pm^{[\kappa]}$, ordered as
\begin{gather*}
  \dots, E^{[\kappa]}_{2},E^{[\kappa]}_{1},
  \SPsi^{[\kappa]}_+,\SPsi^{[\kappa]}_-,
  E^{[\kappa]}_{-1},E^{[\kappa]}_{-2},\dots,
\end{gather*}
and use their appropriate finite subsets to define the $\WW{n[m]}$.
Any such subset is a length-$n$ segment including~$\SPsi_+$
or~$\SPsi_-$ (or both).\footnote{The case where only~$\SPsi_+$, but
  not~$\SPsi_-$ is included is the maximally asymmetric realization
  with $m=0$.  The case where only~$\SPsi_-$ is included,
  but~$\SPsi_+$ is not, is the ``equally asymmetric,'' opposite
  realization with $m=n$.  We do not consider it specially because it
  can be obtained from the $m=0$ realization by the automorphism
  exchanging~$\cE$ and $\cF$.}

Next, the scalar products of $\xi$ with the other vectors
(see~\bref{sec:n[m]}), and hence the OPEs of $Y=(\xi,\dd\varphi)$ with
the other currents are independent of~$k$ or~$n$; this allows us to
\textit{identify the $\YXi$ and $Y$ fields in all the realizations of
  all~$\WW{n}$}.

We now find the dimension-$1$ current in each $\WW{n[m]}$.
\begin{lemma}\label{H-explicit}
  For $0\,{\leq}\, m\,{\leq}\, n$, \ $n\,{\geq}\,2$, the diagonal
  current of $\WW{n[m]}$ is given by
  \begin{multline*}
    \cH_{n[m]}^{(k)}(z)
    = \ell_n(k)\,Y(z)
    +\sum_{i=1}^{n-m-1}\ffrac{n\,{-}\,i\,{-}\,m}{n}\,A^{[k+n]}_i(z)
    + \ffrac{n\,{-}\,m}{n}\,Q^{[k+n]}_+(z)\\*
    - \smash[t]{\sum_{i=-m+1}^{-1}}\ffrac{m\,{+}\,i}{n}\,A^{[k+n]}_i(z)
    - \ffrac{m}{n}\,Q^{[k+n]}_-(z).
  \end{multline*}
\end{lemma}
\begin{proof}
  This is shown by directly finding the centralizer of the relevant
  screening operators in the space of dimension-$1$ operators with
  zero ``momentum'' --- i.e., among descendants of the identity
  operator, and hence necessarily linear combinations of the currents;
  the coefficients are then determined by a straightforward
  calculation, uniquely up to normalization, and the overall
  normalization is fixed by~\eqref{HH-OPE}.
\end{proof}

Proceeding similarly, we establish the existence of a Virasoro algebra
in the centralizer of the screenings:
\begin{lemma}
  For generic $k$, the centralizer of the screenings $E_i$,
  $i=1,\dots,n\,{-}\,1$, and $\SPsi$
  \textup{(}see~\eqref{n[0]-scr}\textup{)} in the space of
  dimension-$2$ operators that are descendants of the identity is
  three-dimensional.  In addition to $\cH\cH(z)$ and $\dd\cH(z)$, it
  is generated by the energy-momentum tensor with central
  charge~\eqref{cn(k)}, explicitly given by
  \begin{multline*}
    \cT_n^{(k)}(z)
    =\fhalf\sum_{i,j\,{\in}\,\{n\,{-}\,1,n-2,\dots,1,{+},{*}\}}
    (\Gamma_{n}(k)^{-1})_{i,j}A^{[k+n]}_iA_j^{[k+n]}(z)\\*
    {}+\sum_{i=1}^{n - 1}
    (n\,{-}\,i)
    \ffrac{(i\,{-}\,1)(k\,{+}\,n\,{-}\,1)\,{-}\,1}{2(k\,{+}\,n)}\,
    \dd A^{[k+n]}_i(z)
    - \ffrac{n}{2}\,\dd A^{[k+n]}_+(z),
  \end{multline*}
  where we write $A_+\equiv Q_+\equiv Q$ and $A_{*}\equiv Y$ for
  notational uniformity and where $\Gamma_n(k)$ is the ``dressed''
  Cartan matrix~\eqref{Gamma}.
  %%   (where $n\ell_n(k)=(n-1)k+n(n-2)$).
\end{lemma}
%% \begin{proof}
%%   Direct calculation.
%% \end{proof}

Lemma~\bref{H-explicit} implies a relation between the $\cH$ currents
in $\WW{n[0]}(k)$ and $\WW{(n-1)[0]}(k+1)$:
\begin{gather}\label{H-recursion}
  n \cH_{n[0]}^{(k)} - (n - 1)\cH_{(n - 1)[0]}^{(k + 1)} =
  (k\,{+}\,n\,{-}\,1)Y + 
  \sum_{i=1}^{n - 1}A_i^{[k+n]} + Q^{[k+n]}.
\end{gather}
Somewhat less obviously, we have similar recursion relations for the
$\cE$ and $\cF$ currents generating $\WW{n[0]}$.

\begin{lemma}\label{basic-recursion}
  Let $\cE_{(n-1)[0]}^{(k+1)}=e^\YXi$, \ $\cH_{(n-1)[0]}^{(k+1)}$, and
  $\cF_{(n-1)[0]}^{(k+1)}=-\polP_{n-1}^{(k+1)}e^{-\YXi}$ be the
  generators of~$\WW{(n-1)[0]}(k\,{+}\,1)$, $n\,{\geq}\,2$.  Then
  \begin{equation*}
    \cF_{n[0]}^{(k)} = \bigl((k\,{+}\,n\,{-}\,1)\dd
    \,{+}\, n \cH_{n[0]}^{(k)}
    \,{-}\, (n \,{-}\, 1)\cH_{n - 1}^{(k+1)}\bigr)
    \cF_{(n - 1)[0]}^{(k+1)},
    \quad n\,{\geq}\,2,
  \end{equation*}
  is in the centralizer of the screenings
  $(E_i)_{i=1,\dots,n\,{-}\,1}$ and~$\SPsi$, i.e., in~$\WW{n[0]}(k)$.
  In other words, $\cF_{n[0]}^{(k)}=-\polP_{n}^{(k)}e^{-\YXi}$ can be
  constructed recursively, via
\begin{equation*}
    \polP_{n}^{(k)}%%(A^{[k+n]}_{n-1},\dots,A^{[k+n]}_{1},Q^{[k+n]}_+)
    =\Bigl(\!(k\,{+}\,n\,{-}\,1)\dd
    \,{+}\,Q^{[k+n]} \,{+}\, \sum_{i=1}^{n-1}A_i^{[k+n]}\Bigr)
    \polP_{n-1}^{(k+1)}
    %%(A^{[k+n]}_{n-2}, \dots, A^{[k+n]}_{1}, Q^{[k+n]}_+)
  \end{equation*}
  with the initial condition %%%~\eqref{rec-bc}.
  $\polP_1^{(k)}(Q)=Q$.
\end{lemma}

In accordance with a remark in~\bref{sec:chores}, this recursion
relies on identification of the fields $A_{n-2}$, \dots, $A_{1}$, $Q$
involved in the construction of $\WW{(n-1)[0]}(k\,{+}\,1)$ with the
corresponding fields among the $A_{n-1},A_{n-2},\dots,A_{1},Q$
involved in the construction of $\WW{n[0]}(k)$.  As a corollary, we note
that the differential polynomial $\polP_{n}$ depends only on
$A_{n-1}$, \dots, $A_{1}$, $Q$.  

\begin{proof}[Proof of~\bref{basic-recursion}]  The second formula in
  \bref{basic-recursion} is equivalent to the first one in view
  of~\eqref{H-recursion}.  We now show that $\cF_{n[0]}^{(k)}$
  constructed via the recursion is in the centralizer of the
  screenings~\eqref{n[0]-scr}.  We assume that this is so for the
  \textit{preceding} generators, in particular,
  \begin{align*}
    [E_i^{[k+n]},\cF_{(n-1)[0]}^{(k+1)}]&=0,\quad i=1,\dots,n\,{-}\,2,
    \\
    [\SPsi^{[k+n]},\cF_{(n-1)[0]}^{(k+1)}]&=0,
  \end{align*}
  or, in terms of OPEs for $i=1,\dots,n\,{-}\,2$,
  \begin{align}\label{Ai-given}
    e^{(a_i,\varphi(z))}\,
    \polP_{n-1}^{(k+1)}(A_{n-2}^{[k+n]},\dots,
    A_{1}^{[k+n]},Q^{[k+n]})(w) &= \ffrac{0}{z\,{-}\,w} \,{+}\, \dots,\\
    \intertext{where the dots denote possible other poles (only the
      vanishing of the indicated pole is essential), and}
    e^{(\psi,\varphi(z))}\,
    \polP_{n-1}^{(k+1)}(A_{n-2}^{[k+n]},\dots,
    A_{1}^{[k+n]},Q^{[k+n]})(w) &= 0\cdot(z\,{-}\,w)^0 \,{+}\, \dots
    \label{Q-given}
  \end{align}
  (that is, the normal-ordered product vanishes).\footnote{For
    $n\,{-}\,1\!=\!1$, with $\polP_{1}(Q)\!=\!Q$, the vanishing
    in~\eqref{Q-given} occurs as follows: the OPE
    \begin{gather*}
      e^{(\psi,\varphi(z))}\,Q(w)=
      -\frac{e^{(\psi,\varphi(z))}}{z\,{-}\,w}
      \,{+}\,\noo{\,e^{(\psi,\varphi)}Q}(w),
    \end{gather*}
    gives zero in the normal-ordered term after the expansion
    $e^{(\psi,\varphi(z))}\!=
    \!e^{(\psi,\varphi(w))}\,{+}\,(z\,{-}\,w)Q(w)$.}  To be precise,
  we should have written $\bs{\varphi}^{[k+n]}$ in the left-hand
  sides. %%% but we ignore this notational purism.
  
  We must show that relations~\eqref{Ai-given} and~\eqref{Q-given}
  imply the OPEs
  \begin{align}\label{Ai-toshow}
    e^{(a_i,\varphi(z))}\,
    \polP_{n}^{(k)}(A_{n-1}^{[k+n]},\dots,
    A_{1}^{[k+n]},Q^{[k+n]})(w) &= \ffrac{0}{z\,{-}\,w} \,{+}\, \dots,
    ~i=1,\dots,n\,{-}\,1,\\
    \intertext{and}
    e^{(\psi,\varphi(z))}\,
    \polP_{n}^{(k)}(A_{n-1}^{[k+n]},\dots,
    A_{1}^{[k+n]},Q^{[k+n]})(w) &= 0\cdot(z\,{-}\,w)^0 \,{+}\, \dots.
    \label{Q-toshow}
  \end{align}
  
  But we readily establish the OPE  
%%   (e.g., from the matrix in~\bref{sec:n[m]}) 
  \begin{equation*}
    A_i(z)\,
    \Bigl(Q%%^{[k+n]}
    (w)
    \,{+}\, \sum_{j=1}^{n - 1}A_j%%^{[k+n]}
    (w)
    \Bigr)= 0,\qquad
    i=1,\dots,n\,{-}\,2,
  \end{equation*}
  and therefore~\eqref{Ai-given} and the recursion for $\polP$
  imply~\eqref{Ai-toshow} for $i\,{=}\,1,\dots,n\,{-}\,2$.  For
  $i=n\,{-}\,1$, we no longer have a similar vanishing OPE, but
  \textit{the only} nonzero OPEs involving $A_{n-1}$ are with itself
  and with~$A_{n-2}$, see~\eqref{rigged}.  The dependence
  of~$\polP_{n}^{(k)}$ on~$A_{n-1}$ is only through the explicit
  occurrences of~$A_{n-1}$ in the recursion relation, and the
  dependence on~$A_{n-2}$ is only through the explicit occurrences
  of~$A_{n-2}$ in the \textit{preceding} recursion relation.  That~is,
  \begin{multline*}
    \polP_n^{(k)}
    =(k\,{+}\,n\,{-}\,1)
    \dd\Bigl(\!(k\,{+}\,n\,{-}\,1)\dd\polP_{n-2}^{(k+2)}
    \,{+}\,\bigl(A_{n-2} \,{+}\,
    \sum_{j=1}^{n-3}A_j \,{+}\, Q\bigr)\polP_{n-2}^{(k+2)}
    \Bigr)\\*
    {}\,{+}\,\bigl(A_{n-1} \,{+}\, A_{n-2} \,{+}\,
    \smash[t]{\sum_{j=1}^{n-3}}
    A_j \,{+}\, Q\bigr)
    \Bigl(\!(k\,{+}\,n\,{-}\,1)\dd\polP_{n-2}^{(k+2)}\qquad\\*
    {}\,{+}\,\bigl(A_{n-2} \,{+}\, \smash[t]{\sum_{j=1}^{n-3}}
    A_j \,{+}\, Q\bigr)\polP_{n-2}^{(k+2)}\Bigr),
  \end{multline*}
  where $\polP_{n-2}^{(k+2)}$ depends only on the fields that have
  zero OPEs with $A_{n-1}$.  This shows that the first-order pole in
  the OPE $e^{(a_{n-1},\varphi(z))}\,\polP_n^{(k)}(w)$
  vanishes; the derivation is elementary.  
%%   , and we only note that it must involve calculating the
%%   \textit{second}-order pole, given~by
%%   \begin{multline*}
%%       {-}(k\,{+}\,n)\frac{\noo{e^{(a_{n-1},\varphi(z))}
%%       \polP_{n-2}^{(k+2)}(w)}}{(z\,{-}\,w)^2}={}\\
%%       {}={-}(k\,{+}\,n)\frac{\noo{e^{(a_{n-1},\varphi)}
%%       \polP_{n-2}^{(k+2)}}\!(w)}{(z\,{-}\,w)^2}
%%       -(k\,{+}\,n)\frac{\noo{e^{(a_{n-1},\varphi)}
%%       A_{n-1}\polP_{n-2}^{(k+2)}}\!(w)}{z\,{-}\,w},
%%   \end{multline*}
%%   which gives the contribution necessary for the first-order pole
%%   cancellation.
  This proves all the equations~\eqref{Ai-toshow}.
  
  It remains to verify~\eqref{Q-toshow}, i.e.,
  $\noo{e^{(\psi,\varphi)}
    \polP_{n}^{(k)}%%%(A_{n-1},\dots,A_{1},Q)
  } =0$.  For this, we note that $Q$ (and hence $e^{(\psi,\varphi)}$)
  has nonvanishing OPEs only with $Q$ and $A_1$; in particular, it
  follows from~\bref{sec:n[m]} that
  \begin{equation*}
    e^{(\psi,\varphi(z))}
    \,\Bigl(Q(w)\,{+}\,\sum_{j=1}^{n-1}A_j(w)\Bigr)=
    (k\,{+}\,n\,{-}\,1)\frac{e^{(\psi,\varphi(w))}}{z\,{-}\,w}.
  \end{equation*}
  The recursion now readily implies that
  $\noo{e^{(\psi,\varphi)}
    \polP_{n}^{(k)}}=(k\,{+}\,n\,{-}\,1)\dd\noo{e^{(\psi,\varphi)}
    \polP_{n-1}^{(k+1)}}=0$ by assumption.
\end{proof}

\begin{example}\label{ex:n=3}
  For $n=3$, it follows from Lemma~\bref{basic-recursion}
  and~\eqref{rec-bc} that
  \begin{multline*}
    \polP_3^{(k)}(A_2,A_1,Q) =
    {A_1}{A_1}{Q}
    \,{+}\, {A_1}{A_2}{Q}
    \,{+}\, 2{A_1}{Q}{Q}
    \,{+}\, A_2 Q{Q}
    \,{+}\, {Q}{Q}{Q}\\*
    \,{+}\, (k\,{+}\,2)\bigl(2 A_1\dd{Q}
    \,{+}\, A_2\dd{Q}  
    \,{+}\, \dd{A_1}{Q}
    \,{+}\, 3\dd{Q}{Q}\bigr)
    \,{+}\, {(k\,{+}\,2)}^2\dd^2{Q}.
  \end{multline*}
  The $n=4$ example is given in Appendix~\ref{sec:W4}.
\end{example}

Generalizing~\bref{basic-recursion} from $\WW{n[0]}$ to $\WW{n[m]}$,
we have
\begin{lemma}\label{basic-recursion-m}
%%   For each $m\,{\geq}\,1$, we have the recursive relations
  Let the currents
  $\;\cE_{(n-1)[m]}^{(k+1)}=\polP^{\dagger(k+1)}_{(n-1)[m]}\,e^\YXi$,
  $\;\cH_{(n-1)[m]}^{(k+1)}$,\, and
  $\;\cF_{(n-1)[m]}^{(k+1)}=(-1)^{m+1}\polP_{(n-1)[m]}^{(k+1)}e^{-\YXi}$
  be the generators of~$\WW{(n\,{-}\,1)[m]}(k\,{+}\,1)$, $n\,{\geq}\, m +
  2$.  Then $\cF_{n[m]}^{(k)}=(-1)^{m+1}\polP_{n[m]}^{(k)}e^{-\YXi}$
  with
  \begin{equation*}
    \polP_{n[m]}^{(k)}
    =\Bigl(\!(k\,{+}\,n\,{-}\,1)\dd
    + Q_+%%^{[k+n]}
    + \!\sum_{i=1}^{n-m-1}\!\!
    A_i%%^{[k+n]}
    \Bigr)
    \polP_{(n-1)[m]}^{(k+1)},\quad
    n\,{\geq}\, m+2,
  \end{equation*}
  and with the initial condition $\polP_{(m+1)[m]}^{(k)}(Q)=Q$ is in
  the centralizer of the screenings
  $(E_i)_{i=1,\dots,n\,{-}\,m\,{-}\,1}$, $\SPsi_+$, $\SPsi_-$,
  $(E_i)_{i=-1,\dots,-m+1}$, i.e., in~$\WW{n[m]}(k)$.
\end{lemma}
We note that this implies that $\polP_{n[m]}^{(k)}$ depends only on
$A_{n-m-1}, \dots, A_{1}, Q_+$; anticipating a similar statement for
$\polP^\dagger$ (that the differential polynomial
$\polP^\dagger_{n[m]}$ depends only on $A_{-1}$, \dots, $A_{-m+1}$,
$Q_-$ for $m\,{\geq}\,2$ (and on $Q_-$ for $m=1$)), we can
express~\eqref{p-a-priori} much more precisely, as
\begin{gather}\label{p-actual}
  \begin{split}
    \cE^{(k)}_{n[m]}(z)&=
    \polP^{\dagger(k)}_{n[m]}(A_{-1}^{[k+n]},\dots,
    A_{-m+1}^{[k+n]},Q_-^{[k+n]})\,e^{\YXi(z)},\\
    \cF^{(k)}_{n[m]}(z)&=
    (-1)^{m+1}\polP^{(k)}_{n[m]}(A_{n-m-1}^{[k+n]},\dots,
    A_{1}^{[k+n]},Q_+^{[k+n]})\,e^{-\YXi(z)}.
  \end{split}
\end{gather}

\begin{proof}[Proof of~\bref{basic-recursion-m}]  The only new element
  compared to the proof of~\bref{basic-recursion} consists in
  verifying the vanishing of the \textit{second}-order pole in the
  operator product involving~$\SPsi_-$,
  \begin{equation*}
    e^{(\psi_-,\varphi(z))}\,
    \polP_{n[m]}^{(k)}(A_{n-m-1}, \dots, A_{1}, Q_+)(w)=
    \ffrac{0}{(z\,{-}\,w)^2}+\dots.
  \end{equation*}
  assuming that this vanishing occurs for $\polP_{(n-1)[m]}^{(k+1)}$.
  We use the notation $[A,B]_n$ for the coefficient at the $n$th-order
  pole in the OPE $A(z)\,B(w)$ and recall the standard relations (see,
  e.g.,~\cite{[Th]})
  \begin{gather*}
    [V,\dd \polP]_2 = [V,\polP]_1+\dd[V,\polP]_2,\\
    [V,[\cA,\polP]_0]_2=[\cA,[V,\polP]_2]_0]
    + \sum_{\ell>0}\,[[V,\cA]_\ell,\polP]_{2-\ell},
  \end{gather*}
  where we take $V=e^{(\psi_-,\varphi)}$, $\cA=Q_+ +
  \!\sum_{i=1}^{n-m-1}\! A_i$, and $\polP=\polP_{(n-1)[m]}^{(k+1)}$.
  The assumption is therefore that $[V,\polP]_2=0$.  It now follows
  from~\bref{sec:n[m]} that $V(z)\,\cA(w)=\frac{[V,\cA]_1}{z\,{-}\,w}$
  with $[V,\cA]_1=-(k\,{+}\,n\,{-}\,1)V$, and therefore
  \begin{multline*}
    [e^{(\psi_-,\varphi)},
    \polP_{n[m]}^{(k)}]_2^{\phantom{y}}
    = [V,(k\,{+}\,n\,{-}\,1)\dd\polP
    + [\cA,\polP]_0]_2^{\phantom{y}}={}\\*
    {}= (k\,{+}\,n\,{-}\,1)[V,\polP]_1 + [[V,\cA]_1,\polP]_{1}=0,
  \end{multline*}
  which is equivalent to the statement that
  $\polP_{n[m]}^{(k)}e^{-\YXi}$, with $\polP_{n[m]}^{(k)}$ given by
  the recursion formula, is in the centralizer of~$\SPsi_-$.
\end{proof}

%% There are other relations among the different $\polP$.
As an immediate consequence of~\bref{basic-recursion}
and~\bref{basic-recursion-m}, we have
\begin{lemma} For
  $1\,{\leq}\, m\,{\leq}\, n\,{-}\,1%%\floor{\frac{n}{2}}
  $,
  \begin{equation*}
    \polP_{n[m]}^{(k)}(A_{n-m-1}, \dots, A_{1}, Q_+)
    = \polP_{(n - 1)[m - 1]}^{(k+1)}(A_{n-m-1}, \dots, A_{1}, Q_+),
  \end{equation*}
  and therefore
  \begin{equation*}
    \polP_{n[m]}^{(k)}
    %%(A^{[k+n]}_{n-m-1}, \dots, A^{[k+n]}_{1}, Q^{[k+n]}_+)
    = \polP_{n - m}^{(k+m)}.
    %%(A^{[k+n]}_{n-m-1}, \dots, A^{[k+n]}_{1}, Q^{[k+n]}_+)
  \end{equation*}
\end{lemma}

Finally, the $\polP^\dagger$ polynomials are also expressed through
the $\polP_{m}$.  This follows by reading the Dynkin diagram
describing the $n[m]$ realization from right to left, which
corresponds to the algebra automorphism interchanging $\cE$ and $\cF$
(and also replacing $m$ with $n-m$).
\begin{lemma}For $m=1,\dots,n\,{-}\,1$,
  \begin{equation*}
    \polP^{\dagger(k)}_{n[m]}
%%     (A^{[k+n]}_{-1},\dots,A^{[k+n]}_{-m},Q^{[k+n]}_-)
    = \polP_{m}^{(k+n-m)}.
%%     (A^{[k+n]}_{-1}, \dots,A^{[k+n]}_{-m},Q^{[k+n]}_-).
  \end{equation*}
\end{lemma}
This shows, in particular, that the $\polP^{\dagger(k)}_{n[m]}$ depend on
the currents as indicated in~\eqref{p-actual}.

\subsection{$\WW{n[m]}$ OPEs}\label{sec:miura}We thus see that all
realizations of all the $\WW{\bullet}$ algebras are determined by a
series of degree-$n$ differential polynomials $\polP_n$ in $n$
variables, $n\,{\geq}\,1$, given by the \textit{normal-ordered} expressions
\begin{multline*}
  \polP_{n}^{(k)}(A_{n-1},\dots,A_{1},Q)(z)={}\\*
  {}=\Bigl(\!(k\,{+}\,n\,{-}\,1)\dd
  +Q(z) + \sum_{i=1}^{n-1}A_i\Bigr)
  \mycirc\Bigl(\!(k\,{+}\,n\,{-}\,1)\dd
  +Q(z) + \sum_{i=1}^{n-2}A_i(z)\Bigr)\mycirc{}\\*
  \mycirc\dots\mycirc
  \Bigl(\!(k\,{+}\,n\,{-}\,1)\dd
  +Q(z) + A_1(z)\Bigr)Q(z).
\end{multline*}
All the $\dd\equiv \frac{\dd}{\dd z}$ operators are applied to the
currents on the right.  Linearly combining the currents as
\begin{alignat*}{2}
  R_0^+&=Q_+,&&\\
  R_i^+&=Q_+ + A_1 +\dots+ A_i,&\quad i&\,{\geq}\,1,\\
  R_0^-&=Q_-,&&\\
  R_i^-&=Q_- + A_{-1} +\dots+ A_{-i},&\quad i&\,{\geq}\,1,
\end{alignat*}
we readily see that their OPEs are indeed those
in~\eqref{RR-ope}--\eqref{RY-ope}, and we obtain
%% For $n\,{-}\,1\,{\geq}\, m\,{\geq}\,1$, the $n[m]$ realization of
%% $\WW{n}$ is thus generated by the currents written in 
Theorem~\bref{thm:free}.
%% \begin{align*}
%%   \cE_{n[m]}^{(k)}&{}={}
%%   \Bigl(\!\bigl((k\,{+}\,n\,{-}\,1)\dd + R_{m-1}^-\bigr)\dots
%%   \bigl((k\,{+}\,n\,{-}\,1)\dd + R_{1}^-\bigr)R_0^-
%%   \Bigr)e^{\YXi}\\
%% \intertext{and}
%%   \cF_{n[m]}^{(k)}&{}={}(-1)^{m+1}
%%   \Bigl(\!\bigl((k\,{+}\,n\,{-}\,1)\dd + R_{n-m-1}^+\bigr)\dots
%%   \bigl((k\,{+}\,n\,{-}\,1)\dd + R_{1}^+\bigr)R_0^+
%%   \Bigr)e^{-\YXi},
%% \end{align*}
%% which generalizes the symmetric $\hSL2$ realization~\eqref{sl2-sym}.
%% The action of the $\dd$ operators is delimited by the brackets.

The other currents in the algebra are to be found from the OPE
$\cE_{n[m]}^{(k)}(z)\,\cF_{n[m]}^{(k)}(w)$.  Calculating it, we first
obtain
\begin{multline*}
  \cE_{n[m]}^{(k)}(z)\,\cF_{n[m]}^{(k)}(w)={}
  (-1)^{m+1}
  e^{\YXi(z)-\YXi(w)}\times{}\\*
  \shoveleft{\qquad
    \Bigl[\bigl((k\,{+}\,n\,{-}\,1)\dd_z
    + R_{m-1}^-(z)+\ffrac{1}{z\,{-}\,w}
    \bigr){}\dots{}}\\*
  \shoveright{\bigl((k\,{+}\,n\,{-}\,1)\dd_z
    + R_{1}^-(z)+\ffrac{1}{z\,{-}\,w}\bigr)
    \bigl(R_0^-(z)+\ffrac{1}{z\,{-}\,w}\bigr)
    \Bigr]\qquad}\\*
  \shoveleft{\qquad\qquad
    \Bigl[\bigl((k\,{+}\,n\,{-}\,1)\dd_w
    + R_{n-m-1}^+(w)-\ffrac{1}{z\,{-}\,w}\bigr){}\dots{}}\\*  
  \bigl((k\,{+}\,n\,{-}\,1)\dd_w + R_{1}^+(w)-\ffrac{1}{z\,{-}\,w}\bigr)
  \bigl(R_0^+(w)-\ffrac{1}{z\,{-}\,w}\bigr)
  \Bigr],
\end{multline*}
where the OPEs with the exponentials have been taken into account and
\textit{it remains to evaluate the OPEs $R_{i}^-(z)R_{j}^+(w)$ between
  the currents in the two (``$z$'' and ``$w$'') brackets}.  The action
of the $\dd$ operators is delimited by the square brackets.  Because
the OPEs $R_{i}^-(z)R_{j}^+(w)$ are independent of $i$ and $j$
(see~\eqref{RR-ope}), we obtain
\begin{multline}
  \cE_{n[m]}^{(k)}(z)\,\cF_{n[m]}^{(k)}(w)={}
  (-1)^{m+1}
  e^{\YXi(z)-\YXi(w)}\\*
  \shoveleft{\qquad
    {}\times\Bigl[\bigl((k\,{+}\,n\,{-}\,1)\dd_z + R_{m-1}^-(z)
    + \ffrac{k\,{+}\,n\,{-}\,1}{(z\,{-}\,w)^2}\nablaplus
    +\ffrac{1}{z\,{-}\,w}
    \bigr){}\dots{}}\\*
  \bigl((k\,{+}\,n\,{-}\,1)\dd_z + R_{1}^-(z)
  + \ffrac{k\,{+}\,n\,{-}\,1}{(z\,{-}\,w)^2}\nablaplus
  + \ffrac{1}{z\,{-}\,w}\bigr)\\*
  \shoveright{\bigl(R_0^-(z)
    + \ffrac{k\,{+}\,n\,{-}\,1}{(z\,{-}\,w)^2}\nablaplus
    + \ffrac{1}{z\,{-}\,w}\bigr)
    \Bigr]\quad\qquad\qquad}\\*
  \shoveleft{\qquad\qquad
    \Bigl[\bigl((k\,{+}\,n\,{-}\,1)\dd_w
    + R_{n-m-1}^+(w)-\ffrac{1}{z\,{-}\,w}\bigr){}\dots{}}\\*
  \bigl((k\,{+}\,n\,{-}\,1)\dd_w + R_{1}^+(w)-\ffrac{1}{z\,{-}\,w}\bigr)
  \bigl(R_0^+(w)-\ffrac{1}{z\,{-}\,w}\bigr)
  \Bigr]\Bigm|_{\nablaplus(1)=0},
\end{multline}
where the right-hand side is normal-ordered and $\nablaplus$ is the
derivation of the ring of differential polynomials such that
\begin{equation*}
  \nablaplus(R_i^+)=1.
\end{equation*}
\textit{After} all $\nablaplus$'s are evaluated in accordance with
this rule, we must set $\nablaplus=0$, which is indicated by the
prescription $\nablaplus(1)=0$.  As noted above, the action of
$\dd_z=\frac{\dd}{\dd z}$ is delimited by the first of the two groups
of factors in brackets (the ``$z$''-factors).  But each $\dd_z$ is
involved only in the combination $\dd_z +
\frac{1}{(z-w)^2}\nablaplus$, where $\nablaplus$ is then applied to
$R_i^+(w)-\frac{1}{z-w}$.  Because
\begin{gather*}
  \ffrac{1}{(z\,{-}\,w)^2}\nablaplus\bigl(R_i^+(w)\bigr)
  =\dd_z \bigl(-\ffrac{1}{z\,{-}\,w}\bigr),
\end{gather*}
it follows that the right-hand side of the OPE
$\cE_{n[m]}^{(k)}(z)\,\cF_{n[m]}^{(k)}(w)$ can be rewritten with all
the $\nablaplus$ dropped and with each $\dd_z$ acting on \textit{all}
factors, either ``$z$'' or ``$w$'', to the right of a given one.  Thus,
\begin{multline}\label{EF-final}
  \cE_{n[m]}^{(k)}(z)\,\cF_{n[m]}^{(k)}(w)={}
  (-1)^{m+1}
  e^{\YXi(z)-\YXi(w)}\\*
  \shoveleft{\quad
    {}\times\bigl((k\,{+}\,n\,{-}\,1)\dd_z + R_{m-1}^-(z)    
    +\ffrac{1}{z\,{-}\,w}
    \bigr){}\dots{}}\\*
  \shoveleft{\qquad\phantom{{}\times{}}\quad
    \bigl((k\,{+}\,n\,{-}\,1)\dd_z + R_{1}^-(z)
    + \ffrac{1}{z\,{-}\,w}\bigr)}
  \bigl((k\,{+}\,n\,{-}\,1)\dd_z
  + R_0^-(z)
  + \ffrac{1}{z\,{-}\,w}\bigr)\\*
  \shoveleft{\qquad\phantom{{}\times{}}\qquad\qquad
    \bigl((k\,{+}\,n\,{-}\,1)\dd_w
    + R_{n-m-1}^+(w)-\ffrac{1}{z\,{-}\,w}\bigr){}\dots{}}\\*
  \bigl((k\,{+}\,n\,{-}\,1)\dd_w + R_{1}^+(w)-\ffrac{1}{z\,{-}\,w}\bigr)
  \bigl(R_0^+(w)-\ffrac{1}{z\,{-}\,w}\bigr).
\end{multline}
In this normal-ordered expression, we expand $e^{\YXi(z)-\YXi(w)}$ up
to the terms $\mathcal{O}\bigl((z\,{-}\,w)^{n-1}\bigr)$ and then
rewrite it as
\begin{gather*}
  \cE^{(k)}_{n[m]}(z)\,\cF^{(k)}_{n[m]}(w)
  = \sum_{j=1}^{n}\frac{\currentU^{(k)}_{n[m],n-j}(w)}{(z\,{-}\,w)^j}.
\end{gather*}
This gives the central term
$\currentU^{(k)}_{n[m],0}=\lambda_{n-1}(n,k)$, the dimension-$1$
current
$\currentU^{(k)}_{n[m],1}(w)=n\lambda_{n-2}(n,k)\cH^{(k)}_{n[m]}(w)$,
the energy-momentum tensor related to $\currentU^{(k)}_{n[m],2}(w)$ as
in~\eqref{XX-gen}, and the other currents
$\currentU^{(k)}_{n[m],i}(w)$, $3\,{\leq}\, i\,{\leq}\, n\,{-}\,1$.

Evaluating the integrals
\begin{gather*}
  \oint dz f(z)
  \bigl(\cE_{n[m]}^{(k)}(z)\,\cF_{n[m]}^{(k)}(w)\bigr)=
  \sum_{j=1}^{n}\mfrac{1}{(j\,{-}\,1)!}\,
  \currentU^{(k)}_{n[m],n-j}(w)\,\dd_w^{j-1}f(w)
\end{gather*}
with suitable test functions, we 
%% express $\opU^{(k)}_{n[m]}$ as (we recall that $\dd_w
%% e^{\YXi(w)}=Y(w)e^{\YXi(w)}$)
arrange the currents $\currentU^{(k)}_{n[m],i}(w)$ into the
order-$(n\,{-}\,1)$ differential operator
\begin{gather*}
  \opU^{(k)}_{n[m]}
  = \sum_{j=1}^{n}\mfrac{1}{(j\,{-}\,1)!}\,
  \currentU^{(k)}_{n[m],n-j}\,\dd^{j-1},
\end{gather*}
similarly to standard cases of the quantum Drinfeld--Sokolov
reduction.  Because of the presence of both derivatives and Cauchy
kernels in the right-hand side of~\eqref{EF-final}, the corresponding
analogue of the Miura mapping is more complicated than in the classic
cases.  We only note that
\begin{align*}
  \opU^{(k)}_{n[m]}
  &=\,
  e^{-\YXi}\opV^{(k)}_{n[m]}(e^{\YXi}{}\cdot{}),
%%   \\
%%   \intertext{where the operator}
%%   \opV^{(k)}_{n[m]}
%%   &=\sum_{j=1}^{n}
%%   \mfrac{1}{(j\,{-}\,1)!}\,
%%   \currentV^{(k)}_{n[m],n-j}\,\dd^{j-1}
\end{align*}
where $\opV^{(k)}_{n[m]}$ is the order-$(n\,{-}\,1)$ differential
operator given by
\begin{multline*}
  (\opV^{(k)}_{n[m]}f)(w)=
  \oint dz\, 
%% (f(z)e^{\YXi(z)})
  f(z)
  \bigl((k\,{+}\,n\,{-}\,1)\ddleft_z - R_{m-1}^-(z)    
  - \ffrac{1}{z\,{-}\,w}
  \bigr){}\dots{}\\*
  \shoveleft{\qquad\phantom{{}\times{}}\quad
    \bigl((k\,{+}\,n\,{-}\,1)\ddleft_z - R_{1}^-(z)
    - \ffrac{1}{z\,{-}\,w}\bigr)
  \bigl((k\,{+}\,n\,{-}\,1)\ddleft_z
  - R_0^-(z)
  - \ffrac{1}{z\,{-}\,w}\bigr)}\\*
  \shoveleft{\qquad\qquad\qquad\quad
    {}\times
%%     e^{-\YXi(w)}
    \bigl((k\,{+}\,n\,{-}\,1)\dd_w
    + R_{n-m-1}^+(w)-\ffrac{1}{z\,{-}\,w}\bigr){}\dots{}}\\*
  \bigl((k\,{+}\,n\,{-}\,1)\dd_w + R_{1}^+(w)-\ffrac{1}{z\,{-}\,w}\bigr)
  \bigl(R_0^+(w)-\ffrac{1}{z\,{-}\,w}\bigr).
\end{multline*}
Derivatives acting on $f$ are written on the right to simplify
comparison with~\eqref{EF-final}.

\begin{example}
  Taking $n=2$ and $m=1$ brings us back to the symmetric realization
  of $\hSL2$.  Equation~\eqref{EF-final} then becomes
  \begin{multline*}
    \cE^{(k)}_{2[1]}(z)\,\cF^{(k)}_{2[1]}(w)={}
    e^{\YXi(z)-\YXi(w)}
    \bigl((k\,{+}\,1)\dd_z + Q_-(z) + \ffrac{1}{z\,{-}\,w}\bigr)
    \bigl(Q_+(w)-\ffrac{1}{z\,{-}\,w}\bigr)\\*
    \qquad\qquad{}=\bigl(1+(z\,{-}\,w)Y
    %%   +\fhalf(z\,{-}\,w)^2\dd Y
    %%   +\fhalf(z\,{-}\,w)^2 YY\bigr)\times{}\\*
    \bigr)
    \Bigl(\ffrac{Q_+-Q_-}{z\,{-}\,w}
    + \ffrac{k}{(z\,{-}\,w)^2}
    \Bigr)\\*
    =\ffrac{k}{(z\,{-}\,w)^2}
    + \ffrac{Q_+ - Q_- + kY}{z\,{-}\,w},
  \end{multline*}
  with $Q_+-Q_-+k Y=2\cH^{(k)}_{2[1]}$.

%% Equation~\eqref{EF-final} for the OPE $\cE_{n[m]}(z)\,\cF_{n[m]}(w)$
%% simplifies for the maximally asymmetric realization,
%% \begin{align*}
%%   \cE_{n}(z)\,\cF_{n}(w)={}&
%%   -e^{\YXi(z)-\YXi(w)}\bigl((k\,{+}\,n\,{-}\,1)\dd_w
%%   + R_{n-1}^+(w)-\ffrac{1}{z\,{-}\,w}\bigr){}\dots{}\\*
%%   &\quad\qquad\bigl((k\,{+}\,n\,{-}\,1)\dd_w +
%%   R_{1}^+(w)-\ffrac{1}{z\,{-}\,w}\bigr)
%%   \bigl(R_0^+(w)-\ffrac{1}{z\,{-}\,w}\bigr).
%% \end{align*}
  Similarly, for the BP algebra ($n=3$) in the maximally asymmetric
  realiation ($m=0$), we evaluate~\eqref{EF-final}
  as%%%\marginlabel{perepisatx cherez $\bs{V}$}
  \begin{align*}
    \cE_{3[0]}(z)\,\cF_{3[0]}(w)={}&
    -e^{\YXi(z)-\YXi(w)}
    \bigl((k\,{+}\,2)\dd_w
    + R^+_2(w) - \ffrac{1}{z\,{-}\,w}\bigr)\\*  
    &\qquad\qquad\qquad\quad
    \bigl((k\,{+}\,2)\dd_w + R^+_1(w) - \ffrac{1}{z\,{-}\,w}\bigr)
    \bigl(R^+_0(w)-\ffrac{1}{z\,{-}\,w}\bigr)\\
    {}={}&
    \bigl(1+(z\,{-}\,w)Y(w)
    +\fhalf(z\,{-}\,w)^2\dd Y(w)
    +\fhalf(z\,{-}\,w)^2 YY(w)
    \bigr)\\*
    &\qquad{}\times\Bigl(
    \ffrac{(2k\,{+}\,3)(k\,{+}\,1)}{(z\,{-}\,w)^3}
    +(k\,{+}\,1)\ffrac{R^+_2(w)\,{+}\,
      R^+_1(w)\,{+}\,R^+_0(w)}{(z\,{-}\,w)^2}\\*
    &\qquad\qquad\qquad{}+\ffrac{(k\,{+}\,2)(\dd R^+_1\,{+}\,2\dd R^+_0)
     \,{+}\,R^+_2R^+_1\,{+}\,R^+_2R^+_0\,{+}\,R^+_1R^+_0}{z\,{-}\,w}
    \Bigr),
  \end{align*}
  with the rest of the calculation totally straightforward.
\end{example}

\section{$\WW{n}$ algebras from $\protect\hSSL{n}{1}$}
\label{sec:from-affine}
The second construction of the $\WW{n}$ algebras is related to the
cosets $\hSSL{n}1/\hSL{n}$, which actually becomes $\WW{n}$ after a
``correction'' by $e^{\pm\sqrt{n}\phi(z)}$, where $\phi$ is an
auxiliary scalar field with the OPE~\eqref{phi-phi}.  This generalizes
the contruction in~\cite{[BFST]} which explicitly shows that after the
appropriate ``correction,'' ${\hSSL21}\!\bigm/\!{\hSL2}$ is $\hSL2$.

Theorem~\bref{thm:coset} is formulated more specifically as follows.
\begin{Thm}\label{thm:coset2}Let the level $k'$ be related to $k$
  by~\eqref{dual-mult}.  Let $e_{1}(z)$, \dots, $e_{n}(z)$ and
  $f_{1}(z)$, \dots, $f_{n}(z)$ be the two $\SL{n}$ $n$-plets of the
  fermionic generators of $\hSSL{n}{1}_{k'}$.  Then the operators
  \begin{align*}
    \cE(z)&=
    \ffrac{1}{(k'\,{+}\,n\,{-}\,1)^{n/2}}\,
    e_{1}(z) e_{2}(z) \dots e_{n}(z)\,
    e^{\sqrt{n}\phi(z)},\\
    \cF(z)&=
    \ffrac{(-1)^{n+1}}{
      (k'\,{+}\,n\,{-}\,1)^{n/2}}\,
    f_{1}(z) f_{2}(z) \dots f_{n}(z)\,
    e^{-\sqrt{n}\phi(z)}
  \end{align*}
  generate the $\WW{n}(k)$ algebra.
\end{Thm}

\subsection{Outline of the proof}\label{sec:outline}

\subsubsection{The second quantum group}\label{2nd-quantum}We begin
with identifying the ``second'' quantum group that commutes with
$\WW{n}$ and with $\qSSL{n}1$ constructed in~\bref{sec:n[m]} (see
footnote~\ref{foot:two} and the preceding text).  For this, we define
\begin{align*}
  S_{n,i}&=\oint e^{-\frac{1}{k+n}(a_i,\varphi)},
  \quad
%%    i=n-m-1,\dots,1,-1,\dots,-m+1,
  i\neq0,\\
  S_{n,0}&=\oint (a\,Q_+ + b\,Q_-) e^{-\frac{1}{k+n}(\psi_+ +
    \psi_-,\varphi)} \quad (a\neq b,~ |a|^2+|b|^2\neq0).
\end{align*}
A direct calculation shows the following lemma.
\begin{lemma}\label{Lemma-S}
  The operators $S_{n,i}$, $n\,{-}\,m\,{-}\,1\,{\geq}\, i\,{\geq}\,
  -m+1$, commute with the $\WW{n[m]}$ algebra and furnish a
  representation of the nilpotent subalgebra of the
  $\Univ_{\tq}\SL{n}$ quantum group with $\tq=e^{i\pi/(k+n)}$.
\end{lemma}

\noindent
For $m=0$ (the maximally asymmetric realization), the operators in the
lemma are $S_{n,i}$, $i=n\,{-}\,1,\dots,1$.  In each of the ``more
symmetric'' realizations, which involve $S_{n,0}$, the arbitrariness
in $a$ and $b$ is due to the possibility of adding a total derivative
to the integrand and of choosing the overall normalization; a certain
pair $(a,b)$ is to be fixed in what follows.

\begin{rem}
  In the notation $\Univ_q\SSL{n}{1}$ used in the previous sections,
  $q$ was part of the general notation for quantum groups.  We now
  have quantum groups with different values of the deformation
  parameter, and we keep the notation $q$ for the deformation
  parameter of~$\SSL{n}{1}$; its \textit{value} is $q=e^{i\pi(k+n)}$,
  expressed through the $k$ paramter in $\WW{n[m]}(k)$ (for all~$m$).
  The quantum deformation of $\SL{n}$ in~\bref{Lemma-S} is with the
  parameter $\tq=e^{\frac{i\pi}{k+n}}$, and we therefore use the
  notation $\Univ_{\tq}\SL{n}$ for the corresponding quantum group.
\end{rem}

For each fixed $n$ and $m$, we next construct vertices that are
highest-weight representations of the $\Univ_{\tq}\SL{n}$ quantum
group generated by the~$S_{n,i}$.  Let
\begin{gather}\label{V-nm}
  V_{n,m}(z)=e^{(v_{n,m}(k),\varphi(z))},
\end{gather}
where
\begin{gather*}
  v_{n,m}(k)=
%%   \ffrac{1}{n(k+n)}\Bigl(
%%   \sum_{i=1}^{n-m-1}(n-i)a_{n-m-i}
%%   +m\psi_+
%%   +m\psi_-
%%   +\sum_{i=1}^{m-1}(m-i)a_{-i}  
%%   \Bigr)
%%   +\ffrac{1}{n}\xi.
  \ffrac{1}{n(k+n)}\Bigl(
  \mathop{\sum_{j=-m+1}^{n-m-1}}_{j\neq0}(m+j)a_j
  +m\psi_+
  +m\psi_-
  \Bigr)
  +\ffrac{1}{n}\xi.
\end{gather*}
The behavior of $V_{n,m}(z)$ under $\Univ_{\tq}\SL{n}$ is described as
follows.  We let $S\ldot V(z)$ denote the ``dressing'' of a vertex
$V(z)$ by a screening operator~$S=\oint du\,s(u)$.  It is given by the
adjoint action of $S$ and can be represented as the integral $\oint
du\,s(u)V(z)$ taken along the contour running below the real axis from
minus infinity to the vicinity of $z$, encompassing~$z$, and returning
to minus infinity above the real axis.
\begin{lemma}
  The vertex operator $V_{n,m}(z)$ satisfies the relations
  \begin{align*}
    S_{n,i}\ldot V_{n,m}(z)&
    =0,\quad n\,{-}\,m\,{-}\,2\,{\geq}\, i\,{\geq}\, {-}m\,{+}\,1
    \kern-50pt\\[-4pt]
    \intertext{and}
    \smash{S_{n,n-m-1}\ldot S_{n,n-m-1}\ldot V_{n,m}(z)}&=0,\\[-4pt]
    \intertext{but}
    \smash{S_{n,n-m-1}\ldot V_{n,m}(z)}&\neq0,
  \end{align*}
  where $S_{n,i}$ are the screenings in Lemma~\ref{Lemma-S}.
  Therefore, $V_{n,m}(z)$ generates the vector representation of
  $\Univ_{\tq}\SL{n}$ \textup{(}the deformation of the vector
  representation of $\SL{n}$\textup{)}.
\end{lemma}
Let $\oC^n_{\tq}$ denote the $n$-dimensional
$\Univ_{\tq}\SL{n}$-representation generated from~$V_{n,m}(z)$.

We next consider the properties of $V_{n,m}(z)$ with respect to the
$\WW{n[m]}$ algebra.  First, because $(\xi,v_{n,m}(k))=0$, \ 
$V_{n,m}(z)$ is \textit{local} with respect to $\cE_{n[m]}$ and
$\cF_{n[m]}$.  Moreover, it is actually primary with respect to
$\cE_{n[m]}$ and $\cF_{n[m]}$.  To formulate this, we consider the
state $\ket{V_{n,m}}$ corresponding to $V_{n,m}$ and introduce the
modes of $\cE_{n[m]}$ and $\cF_{n[m]}$ as
\begin{gather*}
  \cE_{n[m]}(z)=\sum_{\ell\,{\in}\,\oZ+\half\varepsilon_n}
  \cE_{n[m],\ell}\,z^{-\ell-\frac{n}{2}},\qquad
  \cF_{n[m]}(z)=\sum_{\ell\,{\in}\,\oZ+\half\varepsilon_n}
  \cF_{n[m],\ell}\,z^{-\ell-\frac{n}{2}},
\end{gather*}
where $\varepsilon_n=n\,\mathrm{mod}\,2$, and the modes of $\cH$ and
$\cT$ in the standard way, as
\begin{gather*}
  \cH_{n[m]}(z)=\sum_{\ell\,{\in}\,\oZ}\cH_{n[m],\ell}\,z^{-\ell-1},
  \qquad
  \cT_{n[m]}(z)=\sum_{\ell\,{\in}\,\oZ}\cL_{n[m],\ell}\,z^{-\ell-2}.
\end{gather*}
\begin{lemma}
  We have
  \begin{gather}\label{Vnm-kills}
    \begin{alignedat}{2}
      \cE_{n[m],\ell}\,\ket{V_{n,m}}&=0,&\quad
      &\ell\,{\geq}\,{-}\ffrac{n}{2}+1,\\
      \cF_{n[m],\ell}\,\ket{V_{n,m}}&=0,&\quad
      &\ell\,{\geq}\,{-}\ffrac{n}{2}+2,
    \end{alignedat}
  \end{gather}
  and
  \begin{gather}\label{Vnm-eigens}
    \begin{split}
      \cH_{n[m],0}\,\ket{V_{n,m}}&=\ffrac{1}{n}\,\ket{V_{n,m}},\\
      \cL_{n[m],0}\,\ket{V_{n,m}}&=
      %%   -\ffrac{n(n-2)k+(n-3)n^2+1}{2n(k+n)}\,
      \Bigl(1-\ffrac{n}{2}+\ffrac{n^2-1}{2n(k+n)}
      \Bigr)
      \ket{V_{n,m}}.
    \end{split}
  \end{gather}
\end{lemma}
For $n=2$, for example, this gives the standard highest-weight
conditions for $\hSL2$ highest-weight vectors,
$\cE_{\ell\geq0}^{\vphantom{y}}\ket{~}=0$ and
$\cF_{\ell\geq1}^{\vphantom{y}}\ket{~}=0$.

\begin{proof}
  The annihilation conditions in the lemma are shown by directly
  calculating the corresponding OPEs.  We recall that $\cE_{n[m]}(z)$
  and $\cF_{n[m]}(z)$ involve differential polynomials of the
  respective degrees $m$ and $n\,{-}\,m$; a priori, such polynomials
  develop poles of the respective orders $m$ and $n\,{-}\,m$ in the OPE
  with a vertex operator.  But the explicit (factored) form of the
  differential polynomials above readily implies the OPEs
%% The operator products of $\cE_{n[m]}$ and $\cF_{n[m]}(w)$ with
%% $V_{n,m}$ are given by
  \begin{align*}
    \cE_{n[m]}^{(k)}(z)\,V_{n,m}(w)&=0,\\
    \cF_{n[m]}^{(k)}(z)\,V_{n,m}(w)
    &=\ffrac{\widetilde{F}_{n,m}}{z\,{-}\,w},
  \end{align*}  
  which are equivalent to the annihilation conditions in the
  lemma.\footnote{Although we do not need this in what follows, we
    note that $\widetilde{F}_{n,m}(z)=\noo{
%%   \polP^{(k+m+1)}_{n-m-1}(A_{n-m-1}^{[k+m+1]},\dots,
%%     A_1^{[k+m+1]},Q_+^{[k+m+1]})\,e^{-\YXi}
      \cF_{(n-1)[m]}^{(k+1)}\,V_{n,m}}(z)$, the normal-ordered product
    involving the $\cF$ generator of the ``preceding,''
    $\WW{(n-1)[m]}$ algebra, which is realized in the same free-field
    space as discussed in Sec.~\ref{sec:fromqg}.}  The other formulas
  are established immediately.
\end{proof}
%% is therefore $\WW{n[m]}$-primary

For future use, we also note that the length of the momentum of
$V_{n,m}$ is given by
\begin{gather}\label{v-squared}
  (v_{n,m}(k),v_{n,m}(k))
  =\ffrac{n\,{-}\,1}{n(k\,{+}\,n)}.
\end{gather}

\subsubsection{A dual vertex} Similarly to~\eqref{V-nm}, we define the
vertex operator that carries the dual vector representation of
$\Univ_{\tq}\SL{n}$ and is at the same time a twisted
$\WW{n[m]}$-primary.  The results corresponding to the previous two
lemmas are as follows.

For the vertex operator
\begin{gather}\label{V-nm-dual}
  V^*_{n,m}(z)=e^{(v^*_{n,m}(k),\varphi(z))}
\end{gather}
with
\begin{gather*}
  v^*_{n,m}(k)=
  \ffrac{1}{n(k\,{+}\,n)}\Bigl(  
  \mathop{\sum_{j=1}^{n-1}}_{j\neq n-m} j a_{n - m - j}
  + (n\,{-}\,m)\psi_+
  + (n\,{-}\,m)\psi_-
  \Bigr)
  -\ffrac{1}{n}\xi,
\end{gather*}
it follows that
\begin{align*}
  S_{n,i}\ldot V^*_{n,m}(z)&=
  0,\quad n\,{-}\,m\,{-}\,1\,{\geq}\, i\,{\geq}\, {-}m\,{+}\,2
  \kern-50pt\\[-4pt]
  \intertext{and}
  \smash{S_{n,-m+1}\ldot S_{n,-m+1}\ldot V^*_{n,m}(z)}&=0,\\[-4pt]
  \intertext{but}
  \smash{S_{n,-m+1}\ldot V^*_{n,m}(z)}&\neq0,
\end{align*}
showing that $V^*_{n,m}(z)$ generates the dual vector representation
of the quantum group $\Univ_{\tq}\SL{n}$.  We write
$\overset{*}{\oC}{}^n_{\tq}$ for the $n$-dimensional
$\Univ_{\tq}\SL{n}$-module generated from~$V^*_{n,m}(z)$.

Further, $V^*_{n,m}$ is a ``twisted'' $\WW{n[m]}$ primary, namely
\begin{gather}\label{Vnm2-kills}
  \begin{alignedat}{2}
    \cE_{n[m],\ell}\,\ket{V^*_{n,m}}&=0,&\quad
    &\ell\,{\geq}\,{-}\ffrac{n}{2}+2,\\
    \cF_{n[m],\ell}\,\ket{V^*_{n,m}}&=0,&\quad
    &\ell\,{\geq}\,{-}\ffrac{n}{2}+1
  \end{alignedat}
\end{gather}
for the corresponding state.  The last formulas are a reformulation of
the OPEs
\begin{align*}
  \cE_{n[m]}^{(k)}(z)\,V^*_{n,m}(w)
  &=\ffrac{%%\noo{\cE_{n[m+1]}^{(k)}\,V^*_{n,m}
    \widetilde{E}_{n,m}
    %%}
  }{z\,{-}\,w},\\
  \cF_{n[m]}^{(k)}(z)\,V^*_{n,m}(w)
  &=0.
\end{align*}

\subsubsection{$\hSL{n}_{k'}$ vertices and quantum-group
  duality}\label{SLn-vert} We next introduce the $\hSL{n}$ algebra of
the level~${k'}$ related to~$k$ by~\eqref{dual-add} and the
corresponding quantum group~$\Univ_{q'}\SL{n}$.  The vector
representation of $\Univ_{q'}\SL{n}$ can then be ``coupled'' with the
vector representation of $\Univ_{\tq}\SL{n}$, and similarly for the
dual vector representation.  Because we consider generic $k$, we also
have generic $k'$ determined from~\eqref{dual-add}, and hence generic
$\tq=e^{i\pi/(k+n)}$ and~$q'=e^{i\pi/(k'+n)}$ such that $\tq\,q'=-1$.

For the $\hSL{n}$ algebra of level~${k'}$, let $V_{[\lambda_1,k']}(z)$
denote the vertex operator corresponding to the vector representation
weight $\lambda_1$.  We have
\begin{gather*}
  2\rho\cdot\lambda_1=n\,{-}\,1,\qquad
  \lambda_1\cdot\lambda_1=1\,{-}\,\ffrac{1}{n},
\end{gather*}
where $\rho$ is half the sum of positive roots, and hence the
conformal dimension
\begin{align}
  \Delta_{[\lambda_1,k']}
  &= \ffrac{\lambda_1\cdot(\lambda_1+2\rho)}{2(k'\,{+}\,n)}\notag\\
  \intertext{of $V_{[\lambda_1,k']}(z)$ is given by}
  \Delta_{[\lambda_1,k']}
  &= \ffrac{n^2\,{-}\,1}{2n(k'\,{+}\,n)}.
  \label{Delta}
\end{align}
The vertex operator $V_{[\lambda_1,k']}(z)$ has $n$ components, which
make a basis in the vector representation of~$\SL{n}$ (the horizontal
subalgebra in~$\hSL{n}_{k'}$).  We use the
notation~$\oC^n_{\lambda_1}(z)$ for this space, which is the
evaluation representation of $\hSL{n}$ (with the subscript~$\lambda_1$
intended to distinguish it from other copies of $\oC^n$).  In the
Wakimoto bosonization~\cite{[FF]}\pagebreak[3] (also see~\cite{[PRY]} and
references therein), the highest-weight vector in
$\oC_{\lambda_1}^n(z)$ is given by
\begin{gather*}%%%\label{V-Wak}
  V^{(0)}_{[\lambda_1,k']}(z)=
  e^{\frac{1}{\sqrt{k'\,{+}\,n}}\lambda_1\cdot\varphi(z)}.
\end{gather*}
The length squared of its ``momentum'' is
\begin{gather}\label{mom-squared}
  \ffrac{\lambda_1\cdot\lambda_1}{k'\,{+}\,n}
  =\ffrac{n\,{-}\,1}{n(k'\,{+}\,n)}.
\end{gather}

Next, $V_{[\lambda_1,k']}(z)$ carries a representation of the
$\Univ_{q'}\SL{n}$ quantum group, which is the $n$-dim\-ensional
$\Univ_{q'}\SL{n}$-module $\oC^n_{q'}$ given by the quotient of
the Verma module over $n\,{-}\,1$ singular vectors.  The vertex
$V_{[\lambda_1,k']}(z)$ can therefore be represented as
\begin{gather}\label{CCq'}
  \oC_{\lambda_1}^n(z)\tensor\oC^n_{q'}.
\end{gather}

We also introduce the vertex operator
$V^{\phantom{*}}_{[\lambda_{n-1},k']}(z) =V^*_{[\lambda_{1},k']}(z)$
associated with the dual vector representation of~$\SL{n}$, given by
\begin{gather}\label{CC*q'}
  \smash{\overset{*}{\oC}}_{\lambda_{1}}^n(z)
  \tensor\bar{\oC}^n_{q'},
\end{gather}
with the factors dual to the respective factors in~\eqref{CCq'}.  The
dimension and the momentum length squared of the lowest-weight
component coincide with those in~\eqref{Delta}
and~\eqref{mom-squared}.

It is useful to consider the nonlocal algebra $\algA[\SL{n}_{k'}]$ of
vertex operators generated by $\oC_{\lambda_1}^n(z)\tensor\oC^n_{q'}$
and~$\smash{\overset{*}{\oC}}_{\lambda_{1}}^n(z)
\tensor\bar{\oC}^n_{q'}$ in~\eqref{CCq'} and~\eqref{CC*q'}.  It
contains $\Univ\hSL{n}$, and therefore carries the adjoint action
of~$\Univ\hSL{n}$ (in particular, the center acts trivially); in
addition, it carries an action of $\Univ_{q'}\SL{n}$, and for generic
$k'$, \ $\Univ\hSL{n}$ is the space of $\Univ_{q'}\SL{n}$-invariants
in~$\algA[\SL{n}_{k'}]$.

For $\WW{n[m]}$, a similar nonlocal algebra $\algA[\WW{n[m]}(k)]$ is
generated by the vertex operators
\begin{gather*}
  V_{n,m}(z)
  =\oC(z)\tensor\oC^n_{\tq},\qquad
  V^*_{n,m}(z)
  =\overset{*}{\oC}(z)\tensor\bar{\oC}^n_{\tq}.
\end{gather*}
constructed in~\bref{2nd-quantum}.  It also contains $\WW{n[m]}$ (the
$\cE_{n[m]}(z)$ and $\cF_{n[m]}(z)$ currents are identified in the
quantum deformations of the respective products
$\wwedge^{\!n}V_{n,m}(z)$ and $\wwedge^{\!n}V^*_{n,m}(z)$).
%% Therefore, $\algA[\WW{n[m]}(k)]$ carries an action of~$\WW{n[m]}$; 
In addition, $\algA[\WW{n[m]}(k)]$ carries the action
of~$\Univ_{\tq}\SL{n}$, and $\WW{n[m]}$ is the centralizer
of~$\Univ_{\tq}\SL{n}$ for generic~$k$.

\subsubsection{An ``almost local'' subalgebra}
In $\algA[\WW{n[m]}(k)]\tensor\algA[\SL{n}_{k'}]$, with $k$ and $k'$
related by~\eqref{dual-add}, we now identify an ``almost local''
subalgebra $\widetilde{\algL}_{n[m],k}$ by ``coupling'' the special
subspaces in $\algA[\WW{n[m]}(k)]$ and $\algA[\SL{n}_{k'}]$.  That is,
we couple the $\WW{n[m]}$ vertex operators $V_{n,m}(z)$ and
$V^*_{n,m}(z)$ with the $\hSL{n}_{k'}$ vertex operators
$V_{[\lambda_{1},k']}(z)$ and $V^*_{[\lambda_{1},k']}(z)$ as
follows.
%% First, we recall that $V_{n,m}(z)$ and $V^*_{n,m}(z)$ carry
%% $n$-dimensional representations of $\Univ_{\tq}\SL{n}$

Let $R_{\tq}$ and $R_{q'}$ be the $R$-matrices
\begin{align*}
  R_{\tq}:{}\oC^n_{\tq}\tensor\oC^n_{\tq}
  \to{}&\oC^n_{\tq}\tensor\oC^n_{\tq},\\
  R_{q'}:{}\oC^n_{q'}\tensor\oC^n_{q'}\to{}&\oC^n_{q'}
  \tensor\oC^n_{q'}
\end{align*}
for the vector representations of $\Univ_{\tq}\SL{n}$ and
$\Univ_{q'}\SL{n}$ respectively.  We consider the $R$-matrix
\begin{gather*}
  \cR=(R_{\tq})_{13}\tensor(R_{q'})_{24} :
  (\oC^n_{\tq}\tensor\oC^n_{q'})\tensor(\oC^n_{\tq}\tensor\oC^n_{q'})
  \to{}
  (\oC^n_{\tq}\tensor\oC^n_{q'})\tensor(\oC^n_{\tq}\tensor\oC^n_{q'})
\end{gather*}
and recall that~\eqref{dual-add} implies that $\tq q'=-1$.  Then the
tensor product $\oC_{\tq}^{n} {}\tensor{}\oC_{q'}^{n}$ contains a
$1$-dimensional subspace $I_n$ that is invariant under
$\cR^2=\cR_{12}\cR_{21}$ (and in fact, also under~$\cR$).  The
eigenvalue of the thus understood $\cR^2$ operator on this subspace is
$e^{-\frac{2i\pi}{n}}$, which we write~as
\begin{gather}\label{RR}
  \cR^2:I_n\tensor I_n\mapsto e^{-\frac{2i\pi}{n}}I_n\tensor I_n,
\end{gather}

The dual space also contains an invariant $1$-dimensional subspace,
\begin{gather}\label{RR-dual}
  \cR^2:\overset{*}{I}_n\tensor \overset{*}{I}_n
  \mapsto e^{-\frac{2i\pi}{n}}
  \overset{*}{I}_n\tensor \overset{*}{I}_n
\end{gather}
(with $\cR^2$ understood appropriately), and moreover,
\begin{gather}\label{RR-cross}
  \cR^2: I_n\tensor\overset{*}{I}_n
  \mapsto e^{\frac{2i\pi}{n}}
  I_n\tensor\overset{*}{I}_n.
\end{gather}

We next use the embedding $I_n\hookrightarrow\oC_{\tq}^{n}
{}\tensor{}\oC_{q'}^{n}$ to identify an $n$-dimensional subspace in
the tensor product of the vertex operators $V_{n,m}(z)$ and
$V^*_{n,m}(z)$,%%%\marginlabel{make all duals $*$}
\begin{gather}\label{imbed-space}
  \begin{alignedat}{3}
    &\widetilde{\bs{E}}(z)
    \equiv\oC^n(z)\cong
    \oC^n(z)\tensor I_n
    {}\longrightarrow{}
    &&\underbrace{\oC(z){}\tensor{}\oC_{\tq}^{n}}
    &&{}\tensor{}
    \underbrace{\oC_{\lambda_1}^n(z){}\tensor{}\oC_{q'}^{n}}\\[-4pt]
    &
    &&\kern24pt\parallel
    &&\kern42pt\parallel\\[-4pt]
    &
    &&\quad V_{n,m}(z)&&{}\tensor{}~V_{[\lambda_1,k']}(z)
  \end{alignedat}
  \\
  \intertext{and similarly with the dual spaces,}
  \begin{alignedat}{3}
    &\widetilde{\bs{F}}(z)
    \equiv\smash{\overset{*}{\oC}}^n(z)\cong
    \smash{\overset{*}{\oC}}{}^n(z)\tensor I_n
    {}\longrightarrow{}
    &&\underbrace{\smash{\overset{*}{\oC}}{}(z){}\tensor{}
      \bar{\oC}_{\tq}^{n}}
    &&{}\tensor{}
    \underbrace{
      \smash{\smash{\overset{*}{\oC}}{}}_{\lambda_1}^n(z){}\tensor{}
      \bar{\oC}_{q'}^{n}}\\[-4pt]
    &
    &&\kern24pt\parallel
    &&\kern42pt\parallel\\[-4pt]
    &
    &&\quad V^*_{n,m}(z)&&{}\tensor{}~V^*_{[\lambda_1,k']}(z)
  \end{alignedat}
  \label{imbed-space-dual}
\end{gather}
%% Moreover, there exists a trace on the tensor product of
%% $\Univ_q\SL{n}$ and $\Univ_{q^{-1}}\SL{n}$ modules of the
%% same dimension~\cite{[Kass]}.  We use this to take the quantum-group
%% trace of the tensor product of two vertex operators, $V_{n,m}$ and
%% $V_{[\lambda_1,k']}(z)$, and similarly with $V^*_{n,m}$ and
%% $V^*_{[\lambda_1,k']}(z)$:
%% \begin{gather}\label{take-trace}
%%   \begin{split}
%%     \oC(z){}\tensor{}\oC_{q}^{n} {}\tensor{}
%%     \oC{}^n(z){}\tensor{}\oC_{q'}^{n}
%%     &\xrightarrow{\langle\,\underset{2}{{}\cdot{}},
%%       \underset{4}{{}\cdot{}}\,\rangle}
%%     \oC{}^n(z),\\
%%     \overset{*}{\oC}(z){}\tensor{}\bar{\oC}_{q}^{n} {}\tensor{}
%%     \smash{\overset{*}{\oC}}^n(z){}\tensor{}\bar{\oC}_{q'}^{n}
%%     &\xrightarrow{\langle\,\underset{2}{{}\cdot{}},
%%       \underset{4}{{}\cdot{}}\,\rangle}
%%     \smash{\overset{*}{\oC}}^n(z),
%%   \end{split}
%% \end{gather}
In view of the eigenvalues in~\eqref{RR}--\eqref{RR-cross}, which
become the monodromies of the vertex operators constituting $\oC^n(z)$
and $\smash{\overset{*}{\oC}}^n(z)$, these operators are ``almost
local:'' the monodromies of the components of~$\widetilde{\bs{E}}(z)$
and $\widetilde{\bs{F}}(z)$ with respect to each other are
$e^{\pm\frac{2i\pi}{n}}$; the operator products between these
components are
\begin{gather*}
  \widetilde{\bs{E}}_{\alpha}(z)\,
  \widetilde{\bs{E}}_{\beta}(w)
  =(z\,{-}\,w)^{-\frac{1}{n}}\,
  \mathbb{L}^{{+}{+}}_{\alpha\beta}(z,w),\quad
  \widetilde{\bs{F}}_{\alpha}(z)\,
  \widetilde{\bs{F}}_{\beta}(w)
  =(z\,{-}\,w)^{-\frac{1}{n}}\,\mathbb{L}^{{-}{-}}_{\alpha\beta}(z,w),\\
  \widetilde{\bs{E}}_{\alpha}(z)\,
  \widetilde{\bs{F}}_{\beta}(w)
  =(z\,{-}\,w)^{\frac{1}{n}}\,\mathbb{L}^{{+}{-}}_{\alpha\beta}(z,w),
\end{gather*}
where $\mathbb{L}^{{*}{*}}_{\alpha\beta}(z,w)$ are Laurent series in
$(z\,{-}\,w)$.  The ``almost local'' subalgebra
$\widetilde{\algL}_{n[m],k}$ is the algebra generated by
$\widetilde{\bs{E}}$ and $\widetilde{\bs{F}}$.

\subsubsection{A scalar-field ``correction''}We now modify the vertex
operators $\widetilde{\bs{E}}(z)$ and $\widetilde{\bs{F}}(z)$ defined
in~\eqref{imbed-space} and~\eqref{imbed-space-dual} to obtain
fermionic currents.  

We introduce an auxiliary scalar current $\dd f$ with the operator
product
\begin{gather}\label{ff}
  \dd f(z)\,\dd f(w)=\ffrac{1}{(z\,{-}\,w)^2}
\end{gather}
and consider the $1$-dimensional lattice algebra
$\algF_{\frac{1}{\sqrt{n}}}$ generated by
$e^{\pm\frac{1}{\sqrt{n}}f(z)}$.  It then follows that
$\algA[\WW{n[m]}(k)]\tensor\algA[\SL{n}_{k'}]\tensor
\algF_{\frac{1}{\sqrt{n}}}$ contains the local subalgebra generated by
the operators
\begin{align*}
  \bs{E}(z)&=\widetilde{\bs{E}}(z)
  e^{\frac{1}{\sqrt{n}}f(z)},\\
  \bs{F}(z)&=\widetilde{\bs{F}}(z)
  e^{-\frac{1}{\sqrt{n}}f(z)}.
\end{align*}

We note that the dimension of the operators in each of the $n$-plets
given by $\widetilde{\bs{E}}(z)$ and $\widetilde{\bs{F}}(z)$ is the
sum of $\Delta_{[\lambda_1,k']}$ in~\eqref{Delta} and the dimension
$\Delta_{n}(k)$ of $V_{n,m}$ read off from~\eqref{Vnm-eigens}. Because
of~\eqref{dual-add}, we have
\begin{gather}\label{Delta+Delta}
  \Delta_{[\lambda_1,k']}+\Delta_{n}(k)=1-\ffrac{1}{2n}.
\end{gather}
Therefore, the vertex operators $\bs{E}(z)$ and $\bs{F}(z)$ have
dimension~$1$ (given by $1\,{-}\,\frac{1}{2n}$ in~\eqref{Delta+Delta}
plus $\frac{1}{2n}$ for $e^{\pm\frac{1}{\sqrt{n}}f}$).  Similarly, in
a bosonization where (the highest-weight component of) the vertex
operator $\oC_{\lambda_1}^n(z)$ in the right-hand side
of~\eqref{imbed-space} is a pure exponential, the length squared of
its momentum is given by adding~\eqref{v-squared}
and~\eqref{mom-squared},
\begin{gather}\label{length+length}
  \ffrac{n\,{-}\,1}{n(k\,{+}\,n)}
  +\ffrac{n\,{-}\,1}{n(k'\,{+}\,n)}
  =1-\ffrac{1}{n}.
\end{gather}
Therefore, the momentum length squared of $\bs{E}(z)$ and $\bs{F}(z)$
in a bosonized representation is~$1$ (given by $1\,{-}\,\frac{1}{n}$
in~\eqref{length+length} plus $\frac{1}{n}$ for
$e^{\pm\frac{1}{\sqrt{n}}f}$).

To summarize, $\bs{E}(z)$ and $\bs{F}(z)$ are given by $n$ fermionic
dimension-$1$ currents each, which are local with respect to each
other.  $\bs{E}(z)$ is in the vector representation of $\SL{n}$ and
$\bs{F}(z)$ is in the dual representation.  \textit{These $2n$
  dimension-$1$ fermionic currents are the two $n$-plets of the
  $\hSSL{n}{1}$ fermionic currents}, and we thus obtain a
vertex-operator extension of $\WW{n[m]}(k)\tensor\Univ\hSL{n}_{k'}$,
\begin{gather*}
  \xymatrix@=16pt{%
    {}&\widetilde{\algL}_{n[m],k}
    \tensor\algF_{\frac{1}{\sqrt{n}}}
    &{}\\
    \WW{n[m]}(k)\tensor\Univ\hSL{n}_{k'}
%%     \tensor\algF_{\frac{1}{\sqrt{n}}}
    \ar[-1,1]
    &&\Univ\hSSL{n}{1}_{k'}
%%     \tensor\algG_{\sqrt{n}}^{\phantom{y}}
    \ar[-1,-1]
%%     \\
%%     &&\hSSL{n}{1}_{k'}
%%     \ar@{^{(}->}[-1,0]
  }
\end{gather*}
The diagram means that by ``coupling'' vertex operators of $\WW{n}$
and $\hSL{n}$ (and the auxiliary $e^{\pm\frac{1}{\sqrt{n}}f}$), we
obtain an algebra that contains $\hSSL{n}{1}_{k'}$.

This vertex-operator extension ``inverts'' the
coset~\eqref{the-coset}.  This means, in particular, that
$\hSL{n}_{k'}$ in the left-hand part of the above diagram is the
subalgebra of $\hSSL{n}1_{k'}$ in the right, and we actually have the
diagram (dropping the redundant $[m]$)
\begin{gather*}
   \qquad\qquad\xymatrix@=6pt{
    &&\widetilde{\algL}_{n,k}
    \tensor\algF_{\frac{1}{\sqrt{n}}}&&&\\
    \WW{n}(k)\tensor\Univ\hSL{n}_{k'}\tensor\algF_{\frac{1}{\sqrt{n}}}
    \ar[-1,2]
    &&&&\Univ\hSSL{n}{1}_{k'}\tensor\heisenberg_-
    \ar[-1,-2]
    \ar[2,0]
    &&&&&&\\
    &&&&&&\\
    \WW{n}(k)\tensor\Univ\hSL{n}_{k'}\tensor\heisenberg_+
    \ar[-2,0]
    \ar[0,3]
    &&&&\Univ\hSSL{n}{1}_{k'}\tensor\algG_{\sqrt{n}}^{\phantom{y}}
  }
\end{gather*}
The mappings and algebras are here as follows.  The Heisenberg algebra
$\heisenberg_+$ is generated by the current $\dd f(z)$,
see~\eqref{ff}.  By $e^{\pm\frac{1}{\sqrt{n}}f(z)}$, this Heisenberg
algebra is extended to the (nonlocal) lattice algebra
$\algF_{\frac{1}{\sqrt{n}}}$.  The Heisenberg algebra $\heisenberg_-$
is generated by the current $\dd\phi(z)$, see~\eqref{phi-phi}
($\dd\phi(z)$ can be represented as a unique linear combination of the
Cartan currents of $\hSL{n}_{k'}$, the $\cH_{n}$ current, and $\dd f$
that commutes with $\hSSL{n}{1}_{k'}$).  By $e^{\pm\sqrt{n}\phi(z)}$,
it is extended to a one-dimensional lattice vertex-operator algebra
$\algG_{\sqrt{n}}$.  The $\hSL{n}_{k'}$ algebra in the left-hand part
is a subalgebra in $\hSSL{n}{1}_{k'}$.  The algebra
$\Univ\hSSL{n}{1}_{k'} \tensor\algG_{\sqrt{n}}$ contains the elements
$e_{1}(z) \dots e_{n}(z)\tensor e^{\sqrt{n}\phi(z)}$ and $f_{1}(z)
\dots f_{n}(z)\tensor e^{-\sqrt{n}\phi(z)}$ whose images in
$\widetilde{\algL}_{n,k}\tensor\algF_{\!\frac{1}{\sqrt{n}}}$ coincide
with the respective images of the $\cE(z)$ and~$\cF(z)$ currents
of~$\WW{n}$; this describes the horizontal arrow.

\subsection{Lifting the extension}\label{sec:lifting} This
vertex-operator extension can be viewed as resulting from a
Hamiltonian reduction of another vertex-operator extension, namely of
two $\hSL{n}$ algebras with the dual levels related
by~\eqref{dual-add}.  There, we take the $\hSL{n}_k$ and
$\hSL{n}_{k'}$ vertex operators corresponding to the vector
representations, $V_{[\lambda_1,k]}$ and $V_{[\lambda_1,k']}$, and
similarly with the dual operators.  Repeating~\bref{SLn-vert}
for~$\hSL{n}_k$, we view $V_{[\lambda_1,k]}(z)$ as
$\oC^n(z)\tensor\oC^n_{\tq}$; its properties with respect to the
$\Univ_{\tq}\SL{n}$ quantum group are the same as for the $V_{n,m}$
operators of $\WW{n[m]}$.  To distinguish between different $\oC^n$
spaces, we now use the notation $\oC^n_{\lambda_1}(z)$ for
$V_{[\lambda_1,k]}(z)$ and $\oC^n_{\lambda'_1}(z)$ for
$V_{[\lambda_1,k']}(z)$ (which had no prime in the previous
subsection, where we had fewer $\oC^n$ spaces).  {}From~\eqref{RR}, we
again have a $1$-dimensional subspace in
$\oC^n_{\tq}\tensor\oC^n_{q'}$ that is invariant under the
appropriately squared $R$-matrix (and similarly for the dual), and
hence there are $n^2$-dimensional subspaces embedded~as
\begin{align*}
  \oC^n(z)\tensor\smash{\overset{*}{\oC}}^n(z)
  &\longrightarrow
  \oC_{\lambda_1}^n(z)\tensor\oC^n_{\tq}
  \tensor\smash{\overset{*}{\oC}}^n_{\lambda'_1}(z)\tensor\oC^n_{q'},
  \\
  \smash{\overset{*}{\oC}}^n(z)\tensor{\oC}{}^n(z)
  &\longrightarrow
  \smash{\overset{*}{\oC}}^n_{\lambda_1}(z)\tensor\bar{\oC}^n_{\tq}
  \tensor\oC^n_{\lambda'_1}(z)\tensor\bar{\oC}^n_{q'}.
\end{align*}
The thus ``coupled'' vertex operators of $\hSL{n}_k$ and
$\hSL{n}_{k'}$ generate a nonlocal algebra $\widetilde{\algN}\subset
\Univ\hSL{n}_{k}\tensor\Univ\hSL{n}_{k'}$.  Similarly to
$\widetilde{\algL}_{n,k}$, \ $\widetilde{\algN}$ is ``almost local'';
after a free-field ``correction'' as above, we obtain a
vertex-operator algebra $\mathscr{X}_n(k)$ with $2n^2$ fermionic
currents.  Its Hamiltonian reduction then gives $\hSSL{n}{1}$; under
this reduction, the $\oC^{n}_{\lambda_1}(z)$ and
$\smash{\overset{*}{\oC}}^{n}_{\lambda_1}(z)$ vertex operators of
$\hSL{n}$ are reduced, respectively, to the highest-weight vector
$v_+(z)\,{\in}\,\oC(z)$ and the lowest-weight vector
$v_-(z)\,{\in}\,\smash{\overset{*}{\oC}}^{n}(z)$, which are then
identified with the respective $\WW{n}$ vertex operators~\eqref{V-nm}
and~\eqref{V-nm-dual}.  We thus have
\begin{gather*}
  \xymatrix@=10pt{
    {}&
    \widetilde{\algN}\tensor\algF_{\frac{1}{\sqrt{n}}}
    &{}
    \\
    \Univ\hSL{n}_{k}\tensor\Univ\hSL{n}_{k'}%%\oplus\heisenberg
    \ar[-1,1]
    \ar[]+<-30pt,-10pt>;[2,0]+<-30pt,10pt>|{\mathrm{Ham.\ red}.}
    \ar@{=}[]+<22pt,-10pt>;[2,0]+<22pt,10pt>
    &&\mathscr{X}_n(k)
    \ar[-1,-1]
    \ar[2,0]^{\mathrm{Ham.\ red}.}\\
    {}&\widetilde{\algL}_{n,k}\tensor\algF_{\frac{1}{\sqrt{n}}}
%%     \algA[\WW{n}(k)*\hSL{n}_{k'}]
    &{}\\
    \WW{n}(k)\tensor\Univ\hSL{n}_{k'}%%\oplus\heisenberg
    \ar[-1,1]
    &&\hSSL{n}{1}_{k'}%%\oplus\heisenberg
    \ar[-1,-1]
  }
%%   \qquad
%%   \ffrac{1}{k+n} + \ffrac{1}{k'\,{+}\,n} =1,
\end{gather*}
The mappings into the algebra at the top are actually into its local
subalgebras.  For $n=2$, \ $\mathscr{X}_n(k)$ is the exceptional
affine Lie superalgebra $\hD$~\cite{[BFST]}, and the somewhat
mysterious algebras $\mathscr{X}_n(k)$ must therefore give its
``higher'' (W-)analogues.

\section{Conclusions}
For W algebras whose defining set of screenings involves fermionic
screenings, local fields can be explicitly constructed by the method
used in this paper, based on the properties of quantum supergroups.
Here and in~\cite{[FS-D]}, the simplest applications of this technique
have been given; we hope that the method can be applied in much
greater generality.

For the $\WW{n}$ algebras constructed in this paper, it is
interesting to consider integrable representations, which have many
features in common with the integrable representations of~$\hSL2$.
The $\cE(z)$ and $\cF(z)$ currents then satisfy a number of relations
generalizing $\cE(z)^{k+1}=0$ and $\cF(z)^{k+1}=0$ satisfied in the
integrable $\hSL2$ representations.  The corresponding higher
relations are discussed in~\cite{[FJM]}.  Resolutions of the butterfly
type~\cite{[FS-bfly]} are expected to give further insight into the
structure of these representations.  It would also be interesting to
study semi-infinite realizations of such representations, as
in~\cite{[FSt],[semiinf]},
see~\cite{[FJM],[FJMM-Jack],[FJMM-Macdonald]}.\enlargethispage{12pt}

The method of vertex-operator extensions is a powerful tool in
constructing infinite-dimensional algebras and establishing relations
between different algebras; a number of known examples are based on
``duality'' relations of the type of Eq.~\eqref{dual-add} (possibly
with another integer in the right-hand side).  We note that the
construction of ``unifying W~algebras''~\cite{[BEHHH]}, briefly
recalled in the Introduction, involves an even more general relation
between the levels, $\frac{1}{k+n}+\frac{1}{k'+m}=1$, and is worth
being reconsidered from the standpoint of vertex-operator extensions.

\appendix

\section{The lower $\protect\WW{n}$ algebras}\label{app:the-algebras}
%% \addtocounter{subsection}{-1}
The basic operator product for $\WW{n}$ is given by
\begin{multline}\label{XX-gen}
  \cE_n(z)\,\cF_n(w)=\frac{\lambda_{n-1}(n,k)
%     \prod\limits_{i=1}^{n-1}\bigl(i(k\! +\! n\! -\! 1)\! -\! 1\bigr)
  }{
    (z\,{-}\,w)^n}
  {}+\frac{
    n
    \lambda_{n-2}(n,k)
%     \prod\limits_{i=1}^{n-2}\bigl(i(k\! +\! n\! -\! 1)\! -\! 1\bigr)
    \,\cH_n(w)}{(z\,{-}\,w)^{n-1}}\\*
  {}+ \lambda_{n-3}(n,k)
%   \prod\limits_{i=1}^{n-3}\bigl(i(k\,{+}\,n\,{-}\,1) - 1\bigr)
    \frac{
      \frac{n(n\,{-}\,1)}{2}\,\cH_n\cH_n(w)
      \,{+}\,\frac{n((n-2)(k\,{+}\,n\,{-}\,1)-1)}{2}\,\dd\cH_n(w)
      \,{-}\,(k\,{+}\,n)\cT_n(w)
    }{(z\,{-}\,w)^{n-2}}\\*
    \begin{aligned}[t]
      {}+\Bigl(&\lambda_{n-3}(n,k)
      \frac{\cW_{n,3}(w)
        \,{-}\,(k\,{+}\,n)(\half\dd\cT_{n\perp}(w)
        \,{+}\,\frac{1}{\ell_n(k)}\cH_n\cT_{n\perp}(w))
      }{
        (z\,{-}\,w)^{n-3}}
      \\
      &{}+\lambda_{n-2}(n,k)
      \frac{\frac{n}{6 \ell_n(k)^2}\cH_n\cH_n\cH_n(w)
        \,{+}\,\frac{n}{2 \ell_n(k)}\dd\cH_n\,\cH_n(w)
        \,{+}\,\frac{n}{6}\dd^2\cH_n(w)
      }{
        (z\,{-}\,w)^{n-3}}
      \Bigr)      
    \end{aligned}
    \\*
    {}+\mbox{\Large$\dots$},
\end{multline}
where
\begin{gather}\label{lambda-def}
  \lambda_m(n,k)=
  \prod\limits_{i=1}^{m}\bigl(i(k\,{+}\,n\,{-}\,1)\,{-}\,1\bigr),
%%   \qquad  
%%   \lambda^{(2)}_m(n,k)=
%%   \prod\limits_{i=2}^{m}\bigl(i(k\,{+}\,n\,{-}\,1)\,{-}\,1\bigr),
\end{gather}
$\ell_n(k)$ is defined in~\eqref{HH-OPE}, \ $\cT_{n\perp}=\cT_n
-\ffrac{1}{2\ell_n(k)}\cH_n\cH_n$, \ $\cW_{n,3}$ is a Virasoro primary
dimension-$3$ operator with regular OPE with $\cH_n$, and the dots
denote lower-order poles involving operators of
dimensions~\hbox{$\,{\geq}\,4$}.  In what follows, we write the other
operator products for the lower $\WW{n}$ algebras and recall the
realizations of these algebras obtained in accordance
with~\bref{sec:recursion}--\bref{sec:miura}.

\subsection{The $\WW{1}$ algebra} For completeness, we note that
$\WW{1}$ is merely a $\beta\,\gamma$ system\,---\,two \textit{bosons}
with the first-order pole in their OPE.  The $\beta\,\gamma$ system is
well-known to be singled out as the centralizer of one fermionic
screening in the two-boson system and to have two
``realizations''\,---\,bosonizations with either $\beta$ or $\gamma$
given by the exponential times a current.  In accordance with the
recursion established in~\bref{sec:recursion}, this
(current)${}\cdot{}$(exponential) construct is to be encountered in
the symmetric realization of the next algebra, $\WW2$.  This is indeed
the case: the symmetric realization of $\hSL2$ can be obtained by
``rebosonizing'' the Wakimoto representation.

\subsection{The $\hSL2$ algebra}\label{sec:sl2}  For $n=2$, the
$\WW{2}(k)=\hSL2_k$ algebra is singled out from a three-boson system
as the centralizer of either one bosonic and one fermionic or two
fermionic screenings.  This gives the familiar asymmetric and
symmetric three-boson realizations of $\hSL2$.  In the asymmetric
realization (which is $\WW{2[0]}$ in the nomenclature in this paper),
the two currents that generate the algebra are given by
\begin{align*}
  \cE_{2[0]}&=e^\YXi,\\
  \cF_{2[0]}&=-(A_1 Q + Q Q + (k\,{+}\,1) \dd Q)\,e^{-\YXi}
\end{align*}
(see the OPEs at the end of~\bref{sec:n[1]}, where we now set $n=2$).
In the symmetric realization, the $\hSL2$ currents are given by
\begin{equation}\label{sl2-sym}
  \begin{split}
    \cE_{2[1]} &= Q_-\,e^\YXi,\\
    \cF_{2[1]} &= Q_+\,e^{-\YXi}
  \end{split}
\end{equation}
(only the two open dots remain in the rigged Dynkin
diagram~\eqref{rigged}).

These well-known formulas can be easily obtained by directly finding
the centralizer of the corresponding screenings (one fermionic and one
bosonic in the first case and two fermionic in the second case).

\subsection{The BP algebra}\label{sec:BPdef}
By definition~\cite{[Pol],[Ber]}, the BP $W$~algebra $W^{(2)}_3(k) $
is a ``partial'' Hamiltonian reduction of the $\hSL3_k$ affine
algebra, also see~\cite{[dBT]}.  The BP algebra is generated by 
%% (modes of) 
the energy-momentum tensor~$\cT$, two bosonic spin-$\frac{3}{2}$
generators~$\cE$ and~$\cF$, and a scalar current~$\cH$.  In the
operator-product form, the algebra relations are written as (with
$\cX^+=\cE$ and $\cX^-=\cF$ for compactness)
\begin{align*}
  \cE(z)\,\cF(w)={}&\ffrac{(k\,{+}\,1)(2k\,{+}\,3)}{(z\,{-}\,w)^3} +
  \ffrac{3(k\,{+}\,1)\cH(w)}{(z\,{-}\,w)^2}\notag\\*
  &{}+ \ffrac{3\cH\cH - (k\,{+}\,3)\cT
    + \frac{3}{2}(k\,{+}\,1)\dd \cH}{z\,{-}\,w},
%%   \label{XX}
  \\
  \cH(z)\,\cX^\pm(w)={}&\ffrac{\pm \cX^\pm}{z\,{-}\,w},
  \qquad\cT(z)\,\cX^\pm(w)=
  \ffrac{\frac{3}{2}\cX^\pm(w)}{(z\,{-}\,w)^2}
  + \ffrac{\dd \cX^\pm}{z\,{-}\,w},\\
  \cT(z)\,\cH(w)={}& \ffrac{\cH(w)}{(z\,{-}\,w)^2}
  + \ffrac{\dd \cH}{z\,{-}\,w},\\
  \cT(z)\,\cT(w)={}&\ffrac{\half c_3(k)}{(z\,{-}\,w)^4}
  + \ffrac{2 \cT(w)}{(z\,{-}\,w)^2} + \ffrac{\dd \cT}{z\,{-}\,w},\\
  \cH(z)\,\cH(w)={}&\ffrac{\frac{1}{3}(2k\,{+}\,3)}{(z\,{-}\,w)^2}
%%   \label{HH}
\end{align*}
(with regular terms omitted in operator products and with composite
operators ($\cH\cH$) understood to be given by normal-ordered
products).
%% of~\eqref{XX} 

%% The BP algebra proper is the algebra satisfied by the modes of the
%% generators read off from~\eqref{XX}--\eqref{HH}.  The modes are
%% introduced as $\cX^{\pm}(z)=\sum_{m\in\oZ}\cX^{\pm}_m
%% z^{-m-\frac{3}{2}}$, $\cT(z)=\sum_{m\in\oZ}\cL_m z^{-m-2}$, and
%% $\cH(z)=\sum_{m\in\oZ}\cH_m z^{-m-1}$.  The algebra admits a family of
%% automorphisms $\cU_\theta$\,---\,the spectral flow
%% transform\,---\,given by
%% \begin{equation}
%%   \cU_\theta: ~
%%   \begin{aligned}
%%     \cX^{\pm}_m &\mapsto \cX^{\pm}_{m \pm \theta},\\
%%     \cH_m &\mapsto \cH_m + \ffrac{2k\,{+}\,3}{3}\,\theta\delta_{m,0},\\
%%     \cL_m &\mapsto\cL_m + \theta \cH_m
%%     + \ffrac{2k\,{+}\,3}{6}\,\theta^2
%%     \delta_{m,0}
%%   \end{aligned}
%% \end{equation}
%% for $\theta\,{\in}\,\oZ$.

The realizations of $\BP$ are the maximally asymmetric and the ``more
symmetric'' ones (and two more realizations obtained from these two
via the automorphism exchanging $\cE$ and $\cF$).  In the maximally
asymmetric realization,
\begin{gather*}
  \cE^{(k)}_{3[0]}=e^{\YXi},\qquad
  \cF^{(k)}_{3[0]}=-\polP^{(k)}_3(A_2,A_1,Q)e^{-\YXi},
\end{gather*}
with $\polP^{(k)}_3$ given in~\bref{ex:n=3}.  In the ``more
symmetric'' realization, the generators are given by
\begin{gather*}
  \cE^{(k)}_{3[1]}=Q_- e^{\YXi},\qquad
  \cF^{(k)}_{3[1]}=\polP_2^{(k+1)}(A_1,Q_+)e^{-\YXi},
\end{gather*}
where $\polP^{(k)}_2$ is given in~\eqref{rec-bc}.  The OPEs between
the currents involved here are specified in~\bref{sec:n[m]}.

\subsection{The $\WW{4}$ algebra}\label{sec:W4}
\subsubsection{Realizations}The \textit{realizations} of $\WW{4}$,
described in Secs.~\bref{sec:recursion}--\bref{sec:miura}, are as
follows.  In the totally symmetric, $4[2]$, realization, the
generators are given by
\begin{gather*}
  \cE_{4[2]}^{(k)}=\polP_2^{(k+2)}(A_{-1},Q_-)\,e^{\YXi},
  \qquad
  \cF_{4[2]}^{(k)}=-\polP_2^{(k+2)}(A_1,Q_+)\,e^{-\YXi}
\end{gather*}
with $\polP^{(k)}_{2}(A,Q) = A Q + Q Q + (k\,{+}\,1) \dd Q$
(see~\bref{ex:n=2}).  In the $4[1]$ realization, the generators
are
\begin{gather*}
  \cE_{4[1]}^{(k)}=Q_- e^{\YXi},
  \qquad
  \cF_{4[1]}^{(k)}=\polP_3^{(k+1)}(A_2,A_1,Q_+)e^{-\YXi}\!, 
\end{gather*}
where $\polP_3$ is given in~\eqref{ex:n=3}.  Finally, in the maximally
asymmetric realization, the generators are
\begin{gather*}
  \cE_{4[0]}^{(k)}=e^{\YXi},\qquad
  \cF_{4[0]}^{(k)}=-\polP_4^{(k)}(A_3,A_2,A_1,Q)\,e^{-\YXi},
\end{gather*}
where we find from~\bref{basic-recursion} that
\begin{multline*}
  \polP_4^{(k)}(A_3,A_2,A_1,Q)= {A_1}{A_1}{A_1}{Q} +
  2{A_1}{A_1}{A_2}{Q}
  + {A_1}{A_1}{A_3}{Q}\\*
  + 3{A_1}{A_1}{Q}{Q} + {A_1}{A_2}{A_2}{Q} + {A_1}{A_2}{A_3}{Q}
  + 4{A_1}{A_2}{Q}{Q}\\*
  + 2{A_1}{A_3}{Q}{Q} + 3{A_1}{Q}{Q}{Q} + {A_2}{A_2}{Q}{Q}
  + {A_2}{A_3}{Q}{Q}\\*
  + 2{A_2}{Q}{Q}{Q} \shoveright{{}+ {A_3}{Q}{Q}{Q}
    + {Q}{Q}{Q}{Q}}\\*
  \shoveleft{\qquad\quad
    {}+ (k\,{+}\,3)\bigl(3 {A_1}{A_1}\dd{Q} + 4 {A_1}{A_2}\dd{Q} + 2
    {A_1}{A_3}\dd{Q}
    +  {A_1}\dd{A_2}{Q}}\\*
  + 9 {A_1}\dd{Q}{Q} + {A_2}{A_2}\dd{Q} + {A_2}{A_3}\dd{Q}
  + 6 {A_2}\dd{Q}{Q}\quad\\*
  + 3 {A_3}\dd{Q}{Q} + 3 \dd{A_1}{A_1}{Q} + 2 \dd{A_1}{A_2}{Q}
  +  \dd{A_1}{A_3}{Q}\\*
  \shoveright{{}+ 3 \dd{A_1}{Q}{Q} + \dd{A_2}{Q}{Q}
    + 6 \dd{Q}{Q}{Q}\bigr)}\\*
  + {( k\,{+}\,3 ) }^2\bigl(3 \dd{Q}\dd{Q}
  + \dd{A_2}\dd{Q} + \dd^2{A_1}{Q}
  + 3 \dd{A_1}\dd{Q}
  + {A_3}\dd^2{Q}\\*
  + 2 {A_2}\dd^2{Q} + 4 \dd^2{Q}{Q} + 3 {A_1}\dd^2{Q}\bigr)
  + {( k\,{+}\,3 )
  }^3\dd^3{Q}.
\end{multline*}
The remaining $4[3]$ and $4[4]$ realizations follow from $4[1]$ and
$4[0]$ by the automorphism exchanging $\cE$ and $\cF$.

\subsubsection{Operator product expansions}\label{sec:W4OPEs} We now
write the $\WW4$ operator products explicitly.  First,
Eq.~\eqref{XX-gen} has just $4$ pole terms,
\begin{multline*}
%%   \cX^+(z)\,\cX^-(w)
  \cE(z)\,\cF(w)
  =\ffrac{(k\,{+}\,2) (2k\,{+}\,5) (3k\,{+}\,8)}{(z\,{-}\,w)^4}
  +\ffrac{4 (k\,{+}\,2) (2k\,{+}\,5)\cH(w)}{(z\,{-}\,w)^3}\\*
  {} + (k\,{+}\,2)\ffrac{-(k\,{+}\,4)\cT(w)
    + 6 \cH \cH(w)
    + 2  (2k\,{+}\,5)\dd \cH}{(z\,{-}\,w)^2}\\*
  {}+ (k\,{+}\,2)\Bigl(\cW - \thalf (k\,{+}\,4)\dd\cT(w)
    - \ffrac{4 (k\,{+}\,4)}{3k\,{+}\,8}\cT\cH(w)
    + \ffrac{8 (11k\,{+}\,32)}{3 (3k\,{+}\,8)^2}\cH\cH\cH(w)\\*
    {}+ 6 \dd\cH\cH(w)
    + \ffrac{4 (26\,{+}\,17 k\,{+}\,3 k^2)}{3 (3k\,{+}\,8)}\dd^2\cH(w)
  \Bigr)\ffrac{1}{z\,{-}\,w},
\end{multline*}
where $\cT$ is an energy-momentum tensor with the central charge
$c_4(k)$ and $\cW$ is a dimen\-sion-$3$ Virasoro primary whose
operator product with $\cH$ is regular.  In addition to this OPE and
those in~\eqref{HH-OPE}, we have, writing $\cX^+=\cE$ and $\cX^-=\cF$
for compactness, the operator products with $\cW$ given by
\begin{multline*}
  \cW(z)\,\cX^{\pm}(w) = \pm
  \ffrac{2 (k\,{+}\,4) (3k\,{+}\,7) (5k\,{+}\,16)}{(3k\,{+}\,8)^2}\,
  \ffrac{\cX^{\pm}(w)}{(z\,{-}\,w)^3}\\*
  {}+ \ffrac{\pm3\frac{(k+4) (5k+16)}{2 (3k+8)}\,\dd\cX^{\pm}(w)
    - 6\frac{(k+4) (5k+16)}{(3k+8)^2}\,\cH\cX^{\pm}(w)
  }{(z\,{-}\,w)^2}\\*
  {}-\ffrac{k\,{+}\,4}{k\,{+}\,2}\Bigl(
  \ffrac{8(k\,{+}\,3)  }{3k\,{+}\,8} \cH\dd\cX^{\pm}
  + \ffrac{4 (3k^2\,{+}\,15 k\,{+}\,16)} {{(3k\,{+}\,8)}^2}
  \dd\cH\cX^{\pm}
  \mp (k\,{+}\,3)  \dd^2\cX^{\pm}\\*
  \pm \ffrac{2{(k\,{+}\,4)}}
    {3k\,{+}\,8}\cT \cX^{\pm}
  \mp \ffrac{4 (5k\,{+}\,16)}{{(3k\,{+}\,8)}^2}\cH\cH\cX^{\pm}\Bigr)
  \ffrac{1}{z\,{-}\,w},
\end{multline*}
and
\begin{multline*}
  \cW(z)\,\cW(w)=
  \ffrac{2 (k\,{+}\,4) (2k\,{+}\,5) (3k\,{+}\,7)
    (5k\,{+}\,16)}{3k\,{+}\,8}\ffrac{1}{(z\,{-}\,w)^6}\\*
  - \ffrac{(k\,{+}\,4)^2 (5k\,{+}\,16)}{3k\,{+}\,8}\,
  \ffrac{3\Tperp(w)}{(z\,{-}\,w)^4}
  - \ffrac{(k\,{+}\,4)^2 (5k\,{+}\,16)}{2 (3k\,{+}\,8)}\,
  \ffrac{3\dd\Tperp(w)}{(z\,{-}\,w)^3}\\*
  {}+ \Bigl(
  -\ffrac{3(k\,{+}\,4)^2
    (5 k\,{+}\,16) 
    (12k^2\,{+}\,59 k\,{+}\,74)}
  {4(3 k\,{+}\,8) 
    (20k^2\,{+}\,93 k\,{+}\,102)} \dd^2\Tperp(w)\\*
  \shoveright{{}+ \ffrac{8 {(k\,{+}\,4)}^3
      (5 k\,{+}\,16)}{(
      3 k\,{+}\,8) 
      (20k^2\,{+}\,93 k\,{+}\,102)} \Tperp \Tperp(w)
    + 4(k\,{+}\,4)  \Lambda(w)\Bigr)\ffrac{1}{(z\,{-}\,w)^2}}\\*
    {}+\Bigl(
  -\ffrac{(k\,{+}\,4)^2
    (5 k\,{+}\,16) (12k^2\,{+}\,59 k\,{+}\,74)}{6 (3 k\,{+}\,8) 
    (20k^2\,{+}\,93 k\,{+}\,102)}\dd^3\Tperp\\*
  + \ffrac{8(k\,{+}\,4)^3
    (5 k\,{+}\,16)}{(3 k\,{+}\,8) 
    (20k^2\,{+}\,93 k\,{+}\,102)} \dd\Tperp\,\Tperp
  + 2(k\,{+}\,4)  \dd\Lambda\Bigr)\ffrac{1}{z\,{-}\,w},
\end{multline*}
where $\Tperp= \cT - \ffrac{2}{3k\,{+}\,8} \cH \cH$ and
\begin{multline*}
  (k\,{+}\,2)^2\Lambda
  = \cX^+\cX^-
  - \ffrac{k\,{+}\,2}{2} \dd\cW
  - \ffrac{4(k\,{+}\,2)}{3 k\,{+}\,8} \cW \cH
  + \ffrac{3{(k\,{+}\,2)}^2 (k\,{+}\,4)
    (6k^2\,{+}\,33 k\,{+}\,46 )}{2(3 k\,{+}\,8)
    (20k^2\,{+}\,93 k\,{+}\,102)}\dd^2\Tperp\\*
  {} - \ffrac{(k\,{+}\,2) {(k\,{+}\,4)}^2 (11 k\,{+}\,26)}{2
    (3 k\,{+}\,8) (20k^2\,{+}\,93 k\,{+}\,102)} \Tperp \Tperp
  + \ffrac{2(k\,{+}\,2) (k\,{+}\,4)}{3 k\,{+}\,8} \dd(\Tperp \cH)\\*
  {} + \ffrac{8(k\,{+}\,2) (k\,{+}\,4)}{{(3 k\,{+}\,8)}^2} \Tperp \cH \cH
  - \ffrac{8(k\,{+}\,2) (2 k\,{+}\,5)}{3 (3 k\,{+}\,8)} \dd^2\cH \cH
  - \ffrac{2(k\,{+}\,2) (2 k\,{+}\,5)}{3 k\,{+}\,8} \dd \cH\dd\cH\\*
  {} - \ffrac{16(k\,{+}\,2)(2 k\,{+}\,5)}{(3 k\,{+}\,8)^2}\dd\cH\cH\cH
  - \ffrac{32(k\,{+}\,2) (2 k\,{+}\,5)}{3 {(3 k\,{+}\,8)}^3}\cH\cH\cH\cH
  - \ffrac{1}{6}(k\,{+}\,2) (2 k\,{+}\,5) \dd^3\cH
\end{multline*}
is a dimension-$4$ Virasoro primary field.

Next, the operator product of $\cW$ and $\Lambda$ gives rise to a
\textit{single} dimension-$5$ Virasoro primary
$(k\,{+}\,2)^3\cZ=\frac{5}{2}\cX^+\dd\cX^- - \frac{5}{2}\dd\cX^+\cX^-
+ \frac{2}{3k+8}\cH\cX^+\cX^- + \dots$,
\begin{multline*}
  \cW(z)\,\Lambda(w)
  =-\ffrac{12 (k\,{+}\,4) (3 k\,{+}\,5) (3 k\,{+}\,10)
    (4 k\,{+}\,11)}{
    (k\,{+}\,2)(3 k\,{+}\,8) (20 k^2\,{+}\,93 k\,{+}\,102)}\,
  \ffrac{\cW(w)}{(z\,{-}\,w)^4}\\*
  -\ffrac{4(k\,{+}\,4)(36 k^3\,{+}\,279 k^2\,{+}\,695 k\,{+}\,550)}{
    (k\,{+}\,2)
    (3 k\,{+}\,8)
    (20 k^2\,{+}\,93 k\,{+}\,102)}\,
  \ffrac{\dd \cW(w)}{(z\,{-}\,w)^3}\\*
  +(k\,{+}\,4)\Bigl(\!
  (k\,{+}\,4)\cZ(w) + \ffrac{312 (k\,{+}\,4)
    (3 k\,{+}\,5) (3 k\,{+}\,10) 
    (4 k\,{+}\,11)}{
    (k\,{+}\,2)
    (3 k\,{+}\,8) (20k^2\,{+}\,93 k\,{+}\,102)
    (84k^2\,{+}\,349 k\,{+}\,262)} \cW \Tperp(w)\\*
  \shoveright{{}-\ffrac{18
      (3 k\,{+}\,5) (3 k\,{+}\,10) 
      (4 k\,{+}\,11) (4k^2\,{+}\,29 k\,{+}\,62)}{
      (k\,{+}\,2)
      (3 k\,{+}\,8) (20k^2\,{+}\,93 k\,{+}\,102)
      (84k^2\,{+}\,349 k\,{+}\,262)} \dd^2\cW(w)
    \Bigr)\ffrac{1}{(z\,{-}\,w)^2}}\\*
  {} + (k\,{+}\,4)\Bigl(
  \ffrac{2}{5}(k\,{+}\,4)\dd\cZ
  + \ffrac{30(k\,{+}\,4)
    (3 k\,{+}\,5) (3 k\,{+}\,10)}{
    (k\,{+}\,2)^2
    (3 k\,{+}\,8) 
    (84k^2\,{+}\,349 k\,{+}\,262)} \cW\dd\Tperp\\*
  \shoveright{{}-\ffrac{(3 k\,{+}\,5) (3 k\,{+}\,10) 
      (96k^4\,{+}\,1652 k^3\,{+}\,9647 k^2\,{+}\,22746 k\,{+}\,18384)}{
      2(k\,{+}\,2)^2
      (3 k\,{+}\,8) (20k^2\,{+}\,93 k\,{+}\,102)
      (84k^2\,{+}\,349 k\,{+}\,262)}\dd^3\cW\quad}\\*
  {}+ \ffrac{4(k\,{+}\,4) (3 k\,{+}\,5) 
    (3 k\,{+}\,10) (108k^2\,{+}\,523 k\,{+}\,634)}{
    (k\,{+}\,2)^2
    (3 k\,{+}\,8)
    (20k^2\,{+}\,93 k\,{+}\,102) (84k^2\,{+}\,349 k\,{+}\,262)}
  \dd\cW\Tperp
  \Bigr)\ffrac{1}{z\,{-}\,w}.
\end{multline*}

\noindent\textbf{Acknowledgments.} This work was supported in part by
the RFBR Grants 02-02-16946 %%%Fainberg
and 02-01-00930 %%%Batalin
and by the Grant LSS-1578.2003.2; AMS was also supported by the
Foundation for Promotion of Russian Science.

\end{document}